\documentclass[11pt,leqno]{amsart}
\usepackage{epsfig}
\usepackage{amssymb}
\usepackage{amscd}
\usepackage{pstricks}

\usepackage[arrow]{xy}

\newtheorem{theo}{Theorem}[section]
\newtheorem{lemma}[theo]{Lemma}
\newtheorem{defi}[theo]{Definition}
\newtheorem{prop}[theo]{Proposition}
\newtheorem{conj}[theo]{Conjecture}
\newtheorem{cor}[theo]{Corollary}

\theoremstyle{definition}
\newtheorem{remark}[theo]{Remark}

\newtheorem{construction}[theo]{Construction}
\numberwithin{equation}{section}

\def\CP{\mathbb{CP}}
\def\R{\mathbb{R}}
\def\C{\mathbb{C}}
\def\Z{\mathbb{Z}}
\def\PP{{\mathbb P}}
\def\ol{\overline}
\def\O{\mathcal{O}}
\def\cH{\mathcal{H}}

\def\FS{\mathrm{Lag_{vc}}}

\def\bD{{\mathbf D}}
\def\db#1{ \bD^b({#1})}
\def\coh{\operatorname{coh}}
\def\Qcoh{\operatorname{Qcoh}}
\def\lto{\longrightarrow}

\def\Hom{\operatorname{Hom}}
\def\End{\operatorname{End}}
\def\Ext{\operatorname{Ext}}

\def\bdot{{\raise1pt\hbox{$\scriptscriptstyle\bullet$}}}
\def\mod{\mbox{\rm mod--}}
\def\Mod{\mbox{\rm Mod--}}

\def\rk{\operatorname{rk}}
\def\Pic{\operatorname{Pic}}
\def\kk{\mathbb{C}}

\def\D{\mathcal{D}}

\def\T{\mathcal{T}}
\def\Te{\mathcal{E}}
\def\Be{B}

\def\End{\operatorname{End}}
\def\bR{\mathbf{R}}
\def\Bc{\mathfrak{B}}

\def\L{\mathcal{L}}
\def\M{\mathcal{M}}
\def\A{\mathcal{A}}

\def\Ker{\operatorname{Ker}}
\def\Im{\operatorname{Im}}

\title[Mirror symmetry for Del Pezzo surfaces]
{Mirror symmetry for Del Pezzo surfaces:\\
Vanishing cycles and coherent sheaves}

\author{Denis Auroux}
\address{Department of Mathematics, M.I.T., Cambridge MA 02139, USA}
\email{auroux@math.mit.edu}
\author{Ludmil Katzarkov}
\address{Department of Mathematics, University of Miami,
Coral Gables FL 33124, USA}
\email{l.katzarkov@math.miami.edu}
\author{Dmitri Orlov}
\address{
Algebra Section, Steklov Mathematical Institute, Russian Academy
of Sciences, \newline\indent 8 Gubkin str., Moscow 119991, Russia}
\email{orlov@mi.ras.ru}

\thanks{
DA was partially supported by NSF grant DMS-0244844.
LK was partially
supported by NSF grant DMS-9878353 and NSA grant H98230-04-1-0038.
DO was partially supported by  the Weyl Fund, the  grant of the
President of RF in support of young
scientists (No.~MD-2731.2004.1), the Russian Foundation for
Basic Research (No.~05-01-01034), and the Russian Science
Support Foundation.
}

\begin{document}

\begin{abstract}
We study homological mirror symmetry for Del Pezzo surfaces and their mirror
Landau-Ginzburg models. In particular, we show that the derived category
of coherent sheaves on a Del Pezzo surface $X_k$ obtained by blowing up
$\CP^2$ at $k$ points is equivalent to the derived category of vanishing
cycles of a certain elliptic fibration $W_k:M_k\to \C$ with $k+3$ singular fibers,
equipped with a suitable symplectic form. Moreover, we also show that
this mirror correspondence between derived categories can be extended to
noncommutative deformations of $X_k$, and give an explicit correspondence
between the deformation parameters for $X_k$ and the cohomology class
$[B+i\omega]\in H^2(M_k,\C)$. 
\end{abstract}

\maketitle


\section{Introduction}

The phenomenon of mirror symmetry has been studied extensively
in the case of Calabi-Yau manifolds (where it corresponds to a duality
between $N=2$ superconformal sigma models), but also manifests itself in more
general situations. For example, a sigma model whose target space is a Fano
variety is expected to admit a mirror, not necessarily among sigma models, 
but in the more general context of {\it Landau-Ginzburg models}.

For us, a Landau-Ginzburg model is simply a pair $(M,W)$, where $M$ is a
non-compact manifold (carrying a symplectic structure and/or a complex
structure), and $W$ is a complex-valued function on $M$ called
{\it superpotential}. The general philosophy is that, when a
Landau-Ginzburg model $(M,W)$ is mirror to a Fano variety
$X$, the complex (resp.\ symplectic) geometry of $X$ corresponds to the
symplectic (resp.\ complex) geometry of the {\it critical points} of $W$.

We place ourselves in the context of {\it homological mirror symmetry},
where mirror symmetry is interpreted as an equivalence between certain
triangulated categories naturally associated to a mirror pair
\cite{KoICM}. In our case, B-branes on a Fano variety are described by
its derived category of coherent sheaves, and under mirror symmetry they
correspond to the A-branes of a mirror Landau-Ginzburg model.
These A-branes are described by a suitable analogue of the Fukaya category
for a symplectic fibration,
namely the {\it derived category of Lagrangian vanishing cycles}. A rigorous
definition of this category has been proposed by Seidel \cite{Se1} in the case
where the critical points of the superpotential are isolated and
non-degenerate, following ideas of Kontsevich \cite{KoEns} and Hori, Iqbal,
Vafa \cite{HIV}.

Therefore, for a Fano variety $X$ and a mirror Landau-Ginzburg model
$W:M\to\C$, the homological mirror symmetry conjecture can be formulated as
follows:

\begin{conj}\label{conj:hms}
The derived category of Lagrangian vanishing cycles \,$\db{\FS(W)}$ is 
equivalent to the derived category of coherent sheaves \,$\db{\coh(X)}$.
\end{conj}

\begin{remark}
Homological mirror symmetry also predicts another equivalence of derived
categories. Namely, viewing now $X$ as a symplectic manifold and $M$ as a
complex manifold, the derived category of B-branes of the Landau-Ginzburg
model
$W:M\to\C$, which was defined algebraically in \cite{KL,O2} following ideas
of Kontsevich, should be equivalent to the derived Fukaya category of $X$.
This aspect of mirror symmetry will be addressed in a further paper; for
now, we focus exclusively on Conjecture \ref{conj:hms}.
\end{remark}

One of the first examples for which Conjecture \ref{conj:hms} has been
verified is that of $\CP^2$ and its mirror Landau-Ginzburg model
which is the elliptic fibration with three singular fibers determined
by the superpotential $W_0=x+y+1/xy$\, on $(\C^*)^2$ (or rather a fiberwise
compactification of this fibration), see \cite{Se2,AKO}.
Other examples of surfaces for which the derived category of coherent
sheaves has been shown to be equivalent to the derived category of
Lagrangian vanishing cycles of a mirror Landau-Ginzburg model include
weighted projective planes, Hirzebruch surfaces \cite{AKO}, and toric
blow-ups of $\CP^2$ \cite{Ue}. For all these examples, the toric structure
plays a crucial role in determining the geometry of the mirror
Landau-Ginzburg model.

Our goal in this paper is to consider the case of a Del Pezzo
surface $X_K$ obtained by blowing up $\CP^2$ at a set $K$ of $k\le 8$ points
(this is never toric as soon as $k\ge 4$). Our proposal is that a
mirror of $X_K$ can be constructed in the following manner.
Observe that the elliptic fibration with three singular fibers
determined by the superpotential $W_0=x+y+1/xy$\, on $(\C^*)^2$ (i.e., the
mirror of $\CP^2$) admits
a natural compactification to an elliptic fibration
$\ol{W_0}:\ol{M}\to\CP^1$ in which the fiber above infinity consists of
nine rational components (see \S \ref{ss:compactify} for details).
Consider a deformation of $\ol{W_0}$ to another elliptic fibration
$\ol{W_k}:\ol{M}\to\CP^1$, such that $k$ of the $9$ critical points in the
fiber $\ol{W_0}{}^{-1}(\infty)$ are displaced towards finite values of the
superpotential. Let $$M_k=\ol{M}\setminus \ol{W_k}{}^{-1}(\infty),$$ and
denote by $W_k:M_k\to\C$ the restriction of $\ol{W_k}$ to $M_k$.
In the generic case, $W_k$ is an elliptic fibration with $k+3$ nodal fibers,
while $\ol{W_k}{}^{-1}(\infty)$ is a singular
fiber with $9-k$ rational components.
Although we will focus on the Del Pezzo case, this construction also
provides a mirror in some borderline situations. For example, it can be
applied without modification to the case where $\CP^2$ is blown up at
$k=9$ points which lie at the intersection of two
elliptic curves (the fiber $\ol{W_k}{}^{-1}(\infty)$ is then a smooth
elliptic curve).

There are two aspects to the geometry of $M_k$. Viewing $M_k$ as a {\it complex}
manifold (a Zariski open subset of a rational elliptic surface), its complex
structure is closely related to the set of critical values of $W_k$,
which has to be chosen in accordance with a given
symplectic structure on $X_K$. A generic choice of the symplectic structure
on $X_K$ (for which there are no homologically nontrivial Lagrangian
submanifolds) determines a complex structure on $M_k$
for which the $k+3$ critical values of $W_k$ are all distinct (leading to
a very simple category of B-branes). In the opposite situation,
which we will not consider here, if we equip $X_K$ with a
symplectic form for which there are homologically nontrivial Lagrangian
submanifolds, then some of the critical values of $W_k$ become equal, and
the topology of the singular fibers may become more complicated.

The {\it symplectic} geometry of $M_k$ is more important to us. Since
$H^2(M_k,\C)\simeq \C^{k+2}$, the symplectic form $\omega$ on $M_k$, or
rather its complexified variant $B+i\omega$, depends on $k+2$ moduli
parameters. As we will see in \S \ref{s:lagvc}, these parameters completely
determine the derived category of Lagrangian vanishing cycles of $W_k$;
the actual positions of the critical values are of no importance, as
long as the critical points of $W_k$ remain isolated and non-degenerate
(see Lemma \ref{l:exactdeform}). This means that we shall not concern
ourselves with the complex structure on $M_k$; in fact, a compatible
almost-complex structure is sufficient for our purposes, which makes the
problem of deforming the elliptic fibration $\ol{W_0}$ in the prescribed
manner a non-issue.

To summarize, we have:

\begin{construction}
Given a Del Pezzo surface $X_K$ obtained by blowing up $\CP^2$ at $k$
points, the mirror Landau-Ginzburg model is an elliptic
fibration $W_k:M_k\to \C$ with $k+3$ nodal singular fibers, which has
the following properties:

$(i)$ the fibration $W_k$
compactifies to an elliptic fibration $\ol{W_k}$ over $\CP^1$ in which the
fiber above infinity consists of $9-k$ rational components;

$(ii)$ the compactified fibration $\ol{W_k}$ can be obtained
as a deformation of the elliptic fibration $\ol{W_0}:\ol{M}\to\CP^1$ which
compactifies the mirror to $\CP^2$. 

Moreover, the manifold $M_k$ is equipped with a
symplectic form $\omega$ and a B-field $B$, whose cohomology classes are
determined by the set of points $K$ in an explicit manner as discussed
in~\S \ref{s:mainproof}.
\end{construction}

Our main result is the following:

\begin{theo}\label{th:main1}
Given any Del Pezzo surface $X_K$ obtained by blowing up $\CP^2$ at $k$
points, there exists a complexified symplectic form $B+i\omega$ on $M_k$
for which $\db{\coh(X_K)}\cong \db{\FS(W_k)}$.
\end{theo}

The {\it mirror map}, i.e.\ the
relation between the cohomology class $[B+i\omega]\in H^2(M_k,\C)$ and
the positions of the blown up points in $\CP^2$, can be described
explicitly (see Proposition \ref{prop:mirrormap}).

On the other hand, not every choice of $[B+i\omega]\in H^2(M_k,\C)$
yields a category equivalent to the derived category of coherent sheaves
on a Del Pezzo surface. There are two reasons for this. First, certain
specific choices of $[B+i\omega]$ correspond to deformations of the
complex structure of $X_K$ for which the surface contains a
$-2$-curve, which causes the anticanonical class to no longer be ample.
There are many ways in which this can occur, but perhaps the simplest one
corresponds to the case where a same point is blown up twice, i.e.\ we
first blow up $\CP^2$ at $k-1$ generic points and then blow up a point on
one of the exceptional curves. We then say that $X_K$ is obtained from
$\CP^2$ by blowing up $k$ points, two of which are infinitely close, and
call this a ``simple degeneration'' of a Del Pezzo surface.
In this case again we have:

\begin{theo}\label{th:main2}
If $X_K$ is a blowup of $\CP^2$ at $k$ points, two of which are infinitely
close, and a simple degeneration of a Del Pezzo surface, then there exists
a complexified symplectic form $B+i\omega$ on $M_k$
for which $\db{\coh(X_K)}\cong \db{\FS(W_k)}$.
\end{theo}

More importantly, deformations of the symplectic structure on $M_k$ need
not always correspond to deformations of the complex structure on $X_K$
(observe that $H^2(M_k,\C)$ is larger than $H^{1}(X_K,TX_K)$). The
additional deformation parameters on the mirror side can however be
interpreted in terms of {\it noncommutative deformations} of the Del Pezzo
surface $X_K$ (i.e., deformations of the derived category $\db{\coh(X_K)}$).
In this context we have the following theorem, which
generalizes the result obtained in \cite{AKO} for the case of $\CP^2$:

\begin{theo}\label{th:main3}
Given any noncommutative deformation of the Del Pezzo surface $X_K$,
there exists a complexified symplectic form $B+i\omega$ on $M_k$
for which the deformed derived category $\db{\coh(X_{K,\mu})}$ is
equivalent to $\db{\FS(W_k)}$. Conversely, for a generic choice of
$[B+i\omega]\in H^2(M_k,\C)$, the derived category of Lagrangian
vanishing cycles $\db{\FS(W_k)}$ is equivalent to the
derived category of coherent sheaves of a noncommutative deformation of
a Del Pezzo surface.
\end{theo}

The mirror map is again explicit, i.e.\ the parameters which determine
the noncommutative Del Pezzo surface can be read off in a simple manner
from the cohomology class $[B+i\omega]$.

\begin{remark} The key point in the determination of the mirror map is that
the parameters which determine the composition tensors in
$\db{\FS(W_k)}$ can be expressed explicitly in terms of the
cohomology class $[B+i\omega]$
(see \S \ref{ss:m2}). A remarkable feature of these formulas
is that they can be interpreted in terms of {\em theta functions} on a certain
elliptic curve (see \S \ref{ss:theta}). As a consequence, our description
of the mirror map also involves theta functions (see \S \ref{s:mainproof}).
\end{remark}

The rest of the paper is organized as follows. In
\S \ref{s:dbcoh} we describe the bounded derived categories of coherent
sheaves on Del Pezzo surfaces, their simple degenerations, and their
noncommutative deformations. In \S \ref{s:mirrors} we describe the
topology of the elliptic fibration $M_k$ and its vanishing cycles.
In \S \ref{ss:deflagvc} we recall Seidel's definition of the derived
category of Lagrangian vanishing cycles of a symplectic fibration, and
in the rest of \S \ref{s:lagvc} we determine $\db{\FS(W_k)}$. Finally
in \S \ref{s:mainproof} we compare the two viewpoints, describe the mirror
map, and prove the main theorems.

\bigskip

\noindent
{\bf  Acknowledgements:}
We are thankful to A. Kapustin, T. Pantev, P. Seidel for many helpful discussions.


\section{Derived categories of coherent sheaves on blowups of
$\CP^2$} \label{s:dbcoh}

The purpose of this section is to give a description of the bounded
derived categories of coherent sheaves on Del Pezzo surfaces, their
simple degenerations, and
their noncommutative deformations. 
We always work over the field of complex numbers $\kk.$

\subsection{Del Pezzo surfaces and blowups of the projective plane at
distinct points}\label{ss:delpezzo}

\begin{defi} A smooth projective surface $S$ is called
a Del Pezzo surface if the anticanonical sheaf $\O_S(-K_S)$ is
ample (i.e., a Del Pezzo surface is a Fano variety of dimension 2).
\end{defi}

The Kodaira vanishing theorem and Serre duality give us immediately
that for any Del Pezzo surface
\begin{align*}
&H^1(S, \O(-mK_S))=0
\qquad
\text{for all } \ m\in \Z,\\
&H^2(S, \O(-mK_S))=0
\qquad
\text{for all } \ m\ge 0,\\
&H^2(S, \O(-mK_S))=H^0(S, \O((m+1)K_S))
\qquad
\text{for all } \ m\in \Z.
\end{align*}
In particular, we obtain that $H^1(S, \O_S)\cong H^2(S, \O_S)=0,$ and
$H^0(S, \O(mK_S))=0$ for all $m>0.$
By the Castelnuovo-Enriques  criterion any Del Pezzo surface is rational.

Let $S$ be a Del Pezzo surface. The integer $K^2_S$  is called the {\em degree} of $S$
and will be denoted by~$d$.
The Noether formula gives a relation between the degree and the rank of the Picard group of a
Del Pezzo surface:
$
d=K^2_S = 10- \rk \Pic S \le 9.
$

We  can also introduce another integer number which is called the {\em index}
of $S$. This is the maximal $r>0$ such that $\O(-K_S)=\O(rH)$
for some divisor $H$. The inequality  $d\le 9$ implies that $r\le 3.$

Now recall the classification of   Del Pezzo surfaces.

If $r=3,$ then  $S\cong\PP^2$ is the projective plane and $d=9.$
If $r=2,$ then  $S\cong\PP^1\times \PP^1$ is the quadric and $d=8.$
The other Del Pezzo surfaces are not minimal and can be obtained
by blowing up the projective plane $\PP^2.$ More precisely, if $S$ is a Del Pezzo
surface of index $r=1,$ then it has degree $1\le d\le 8$ and
$S$ is a blowup of the projective plane $\PP^2$ at
$k=9-d$ distinct points. The ampleness of the anticanonical class
requires that in this set no three points lie on a line, and no six
points lie on a conic; moreover, if $k=8$ the eight points are not allowed
to lie on an irreducible cubic which has a double point at one of these
points.
Conversely, any surface which is
a blowup of the projective plane at a set of $k\le 8$ different points
satisfying these constraints is
a Del Pezzo surface of degree $d=9-k.$
All these facts are well-known and can be found in any textbook on
surfaces (see e.g.\ \cite{GH}).

\medskip

Denote by $\db{\coh(S)}$ the bounded derived category of coherent sheaves on $S$.
It is known that the bounded derived category of coherent sheaves on any
Del Pezzo surface has a full exceptional collection, which makes it possible
to establish
an equivalence between the category $\db{\coh(S)}$ and the bounded derived category
of finitely generated modules over the algebra of the exceptional collection (\cite{Or}, see
also \cite{KuOr}). This is a particular case of a more general statement about
derived categories of blowups.

First, recall the notion of exceptional collection.
\begin{defi}
An object $E$ of a $\kk$-linear triangulated category $\D$ is
said to be {\sf exceptional} if
$\Hom(E, E[k])=0$ for all $k\ne 0,$ and $\Hom(E, E)=\kk$.
An ordered set of exceptional objects $\sigma=\left(E_0,\ldots E_n\right)$
is called an {\sf exceptional collection} if
$\Hom(E_j, E_i[k])=0$ for $j>i$ and all $k.$
The exceptional collection $\sigma$ is said to be {\sf strong} if  it satisfies
the additional condition $\Hom(E_j, E_i[k])=0$ for all $i,j$ and for $k\ne 0.$
\end{defi}
\begin{defi} An exceptional collection
$\left(E_0,\ldots,E_n\right)$ in a category $\D$ is called {\sf full}
if it generates the category
$\D,$ i.e. the minimal triangulated subcategory of $\D$
containing all objects $E_i$ coincides
with $\D.$ In this case we say that $\D$ has a semiorthogonal decomposition
of the form
$$
\D=\left\langle E_0,\ldots,E_n\right\rangle.
$$
\end{defi}

The most studied example of an exceptional collection is
the sequence of invertible sheaves
$\langle \O_{\PP^n},\dots,\O_{\PP^n}(n)\rangle$ on the projective space $\PP^n$ (\cite{Be}).
In particular, this exceptional collection on the projective plane $\PP^2$
has length $3$.

\begin{defi} \label{def:algexc}
The  algebra  of a strong exceptional collection
$\sigma=\left(E_0,\ldots,E_n\right)$
is the algebra of endomorphisms $\Be(\sigma)=\End(\Te)$ of the object
$\Te=\mathop\oplus\limits_{i=0}^{n}E_i.$
\end{defi}

Assume that the triangulated category $\D$ has a full strong
exceptional collection $\left(E_0,\ldots, E_n\right)$ and $\Be$ is
the corresponding algebra. Denote by $\mod\Be$ the category of
finitely generated right modules over $\Be.$ There is a theorem
according to which if $\D$ is an {\it enhanced triangulated category}
in the sense of Bondal and Kapranov \cite{BK1}, then
it is equivalent to the bounded
derived category $\db{\mod\Be}$. This equivalence
is given by the functor $\bR\Hom(\Te, -)$  (see \cite{BK1}).

For example, if $\D\cong \db{\coh(X)}$ is the bounded derived category of coherent
sheaves on a projective variety $X,$ then it is enhanced.
Actually, the category of quasi-coherent sheaves $\Qcoh$
has enough injectives, and
$\db{\coh(X)}$ is equivalent to
 the full subcategory $\bD^b_{\coh}(\Qcoh(X))\subset
\db{\Qcoh(X)}$ whose objects are complexes with cohomologies in
$\coh(X).$

Assume that $X$ is smooth and $(E_0,\ldots, E_n)$ is a strong exceptional collection on $X.$
The object $\Te=\bigoplus_{i=0}^n E_i$ defines the derived functor
$$
\bR\Hom(\Te, -):\bD^{+}(\Qcoh(X))\lto \bD^{+}(\Mod\Be),
$$
where $\Mod\Be$ is the category of all right modules over $\Be.$
Moreover, the functor $\bR\Hom(\Te, -)$ sends objects of
$\bD^b_{\coh}(\Qcoh(X))$ to objects of the subcategory
$\bD^b_{\operatorname{mod}}(\Mod\Be),$ which is also equivalent to
$\db{\mod\Be}.$ This gives us a functor
$$
\bR\Hom(\Te, -):\db{\coh(X)}\lto \db{\mod\Be}.
$$
The objects $E_i$ for $i=0,\ldots, n$ are mapped to the projective
modules $P_i=\Hom(\Te, E_i).$ Moreover,
$\Be=\mathop\bigoplus_{i=0}^{n} P_i.$ The algebra $\Be$ has $n+1$
primitive idempotents $e_i,\; i=0,\ldots,n$ such that
$1_\Be=e_0+\cdots+e_{n}$ and $e_i e_j=0$ if $i\ne j.$ The right
projective modules $P_i$ coincide with $e_i \Be.$ The morphisms
between them can be easily described since
$$
\Hom(P_i, P_j)=\Hom(e_i \Be, e_j \Be)\cong e_j \Be e_i\cong
\Hom(E_i, E_j).
$$
This yields an equivalence between the triangulated subcategory of
$\db{\coh(X)}$ generated by the collection $\langle
E_0,\ldots,E_n\rangle$ and the derived category $\db{\mod B}.$
Here we use the fact that the algebra $\Be$ has  a finite global
dimension and any right (and left) module $M$ has a finite
projective resolution consisting of the projective modules $P_i$
with $i=0,\dots, n.$
Finally, if the collection $(E_0,\dots, E_n)$ is full, then we
obtain an equivalence between $\db{\coh(X)}$ and $\db{\mod B}.$

Sometimes it is useful to represent the algebra $\Be$ as a
category $\Bc$ which has $n+1$ objects, say $v_0,\ldots,v_{n},$ and
morphisms  defined by the rule
$
\Hom(v_i, v_j)\cong\Hom(E_i, E_j)
$
with the natural composition law. Thus
$\Be=\mathop\bigoplus\limits_{0\le i,j\le n}\Hom(v_i, v_j).$

\begin{theo}{\rm \cite{Or, KuOr}}\label{fecBu}
Let $\pi: X_K\to \PP^2$ be a blowup of the projective plane
$\PP^2$ at a set $K=\{p_1,\dots, p_k\}$ of any $k$
distinct points, and let $l_1,\dots, l_k$ be the exceptional curves of
the blowup. Let $(F_0, F_1, F_2)$ be a full strong
exceptional collection of vector bundles on $\PP^2.$ Then the
sequence
\begin{equation}\label{excol}
\left(\pi^* F_0,\pi^* F_1, \pi^* F_2, \O_{l_1},\dots,\O_{l_k}\right),
\end{equation}
where the $\O_{l_i}$ are the structure sheaves of the exceptional
$-1$-curves $l_i$, is a full strong exceptional collection on
$X_K$. Moreover, the sheaves $\O_{l_i}$ and $\O_{l_j}$ are
mutually orthogonal for all $i\ne j.$

In particular, there is an equivalence
\begin{equation}\label{equiv:del}
\db{\coh(X_K)}\cong\db{\mod B_K},
\end{equation}
where $B_K$ is the algebra of
homomorphisms of the exceptional collection (\ref{excol}).
\end{theo}
There are no restrictions on the  set of points $K=\{p_1,\dots,p_k\}$ in
this theorem and, in particular, we do not need to assume that $X_K$ is a
Del Pezzo surface.

We can easily describe the space of morphisms from $\pi^* F_i$ to the
sheaf $\O_{l_j}$, since it is naturally identified with the space that
is dual to the fiber of the vector bundle $F_i$ at the point
$p_j\in \PP^2,$ i.e.
$$
\Hom_{X_K}(\pi^* F_i, \O_{l_j})\cong\Hom_{\PP^2}(F_i, \O_{p_j}).
$$

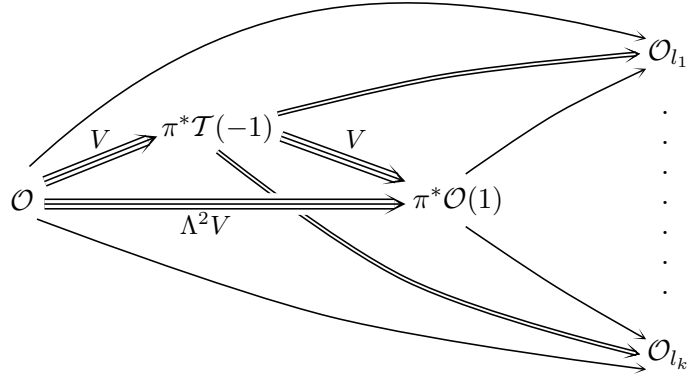
\begin{figure}[t]
\setlength{\unitlength}{1mm}
\begin{picture}(85,52)(30,-22)
\psset{unit=\unitlength,arrowsize=4pt 1,arrowlength=1,arrowinset=0.75,linewidth=0.6pt}
\put(29,0.5){\makebox(0,0)[cc]{$\O$}}
\put(55,10){\makebox(0,0)[cc]{$\pi^*\T(-1)$}}
\put(87,0.5){\makebox(0,0)[cc]{$\pi^*\O(1)$}}
\put(115,20){\makebox(0,0)[cc]{$\O_{l_1}$}}
\multiput(114,12)(0,-4){7}{.}
\put(115,-20){\makebox(0,0)[cc]{$\O_{l_k}$}}
\psbezier{->}(30,5)(57,30)(80,30)(112,22)
\pscurve[doubleline=true]{->}(55,7)(80,-10)(111.2,-20)
\psline[doubleline=true,doublesep=5\pslinewidth,border=2pt]{-}(32,0)(78.5,0)
\psline[linewidth=\pslinewidth,arrowsize=4pt 7]{->}(32,0)(80,0)
\psline[doubleline=true,doublesep=5\pslinewidth]{-}(32,3)(45.75,8.5)
\psline[linewidth=\pslinewidth,arrowsize=4pt 7]{->}(32,3)(47,9)
\psline[doubleline=true,doublesep=5\pslinewidth]{-}(63.5,9)(78.8,3.436)
\psline[linewidth=\pslinewidth,arrowsize=4pt 7]{->}(63.5,9)(80,3)
\pscurve[doubleline=true]{->}(63,12)(85,17)(111.2,20)
\pscurve{->}(88,4)(100,12)(112,18)
\pscurve{->}(88,-3)(100,-11)(112,-18)
\pscurve{->}(31,-2)(70,-15)(112,-22)
\put(50,-4){\small $\Lambda^2 V$}
\put(38,7.5){\small $V$}
\put(72,7.5){\small $V$}
\end{picture}
\caption{The quiver $\Bc_K$ for a blowup of $\PP^2$ at $k$ distinct points.}
\label{quiver}
\end{figure}

There are various standard exceptional collections on the projective plane.
One of them is the collection of line bundles $(\O,\O(1),\O(2)),$
another is the collection $(\O, \T_{\PP^2}(-1), \O(1)),$
where $\T_{\PP^2}$ is the tangent bundle on $\PP^2.$
The latter choice is the most convenient for us. It is easy to see that
$$
\Hom(\O, \T_{\PP^2}(-1))\cong\Hom(\T_{\PP^2}(-1), \O(1))\cong V
\qquad
\text{and}
\qquad
\Hom(\O, \O(1))\cong \Lambda^2 V\cong V^*,
$$
where $V$ is the 3-dimensional vector space whose projectivization $\PP(V)$
is the given projective plane~$\PP^2.$

Let us consider the blowup $X_K$ of the projective plane $\PP(V)$ at a set
$K=\{p_1,\dots,p_k\}$ of $k$ distinct points, and
the exceptional collection
\begin{equation}\label{sigma}
\sigma=(\O_{X_K}, \pi^*\T_{\PP^2}(-1), \pi^*\O_{\PP^2}(1), \O_{l_1},\dots, \O_{l_k}).
\end{equation}
Let $\Bc_K(\sigma)$ be the category of homomorphisms
of this exceptional collection (see Figure \ref{quiver}).
Then the surface $X_K$ can be recovered from the category $\Bc_K(\sigma)$
by means of the following procedure.

Denote by $S_j$ the 2-dimensional space of homomorphisms
from $\pi^*\T_{\PP^2}(-1)$ to $\O_{l_j}$ and denote by $U_j$ the 1-dimensional
space of homomorphisms from $\O(1)$ to $\O_{l_j}.$ The composition law in the category
$\Bc_K(\sigma)$ gives a map from $U_j\otimes V$ to $S_j.$
The kernel of this map is a 1-dimensional subspace $V_j\subset V$, which defines
a point $p_j\in\PP(V)$. In this way, we can determine all the points
$p_1,\dots, p_k\in \PP(V)$
and completely recover the surface $X_K$ starting from the category 
$\Bc_K(\sigma).$

\begin{remark}
Exceptional objects and exceptional collections
on Del Pezzo surfaces are well-studied objects.
First, any exceptional object of the derived category is isomorphic to a sheaf up to translation.
Second, any exceptional sheaf can be included in a full exceptional collection.
Third, any full exceptional collection can be obtained from a given one by
a sequence of natural operations on exceptional collections called {\it
mutations}. All these facts can be found in the paper \cite{KuOr}.
\end{remark}

\subsection{Simple degenerations of Del Pezzo surfaces}\label{ss:simpdeg}

We now look at some simple degenerations of the situation considered above,
namely when two points, for example $p_1$ and $p_2,$ converge to
each other and finally coincide. More precisely, 
this means that we first blow up a point $p$ and after that we blow up
some point $p'$ on the $-1$-curve which is the pre-image of $p$ under the
first blowup. This operation is sometimes called a blowup at two
``infinitely close'' points; more precisely, it corresponds to blowing up
a subscheme of length $2$ supported at $p$.
In this case, the pre-image $\pi^{-1}(p)$ consists of two rational curves meeting
at one point. One of them is a $-1$-curve which we denote by $l'$, and the other
is a $-2$-curve which we denote by $l$. The curve $l$ is the proper
transform of the exceptional curve of the first blow up performed at the
point $p\in \PP^2.$

In this paragraph, we consider the situation where the surface $X_K$ is
the blowup of the projective plane $\PP^2$ 
at a subscheme $K$ which is supported at a set of
$k-1$ points $\{ p, p_3,\dots, p_k\}$ and has length 2 at the point $p$.
In this case the surface $X_K$ is not a Del Pezzo surface, because it
possesses a $-2$-curve $l$.
However, it follows from general results about blowups that $\db{\coh(X_K)}$
still possesses a full exceptional collection \cite{Or}.

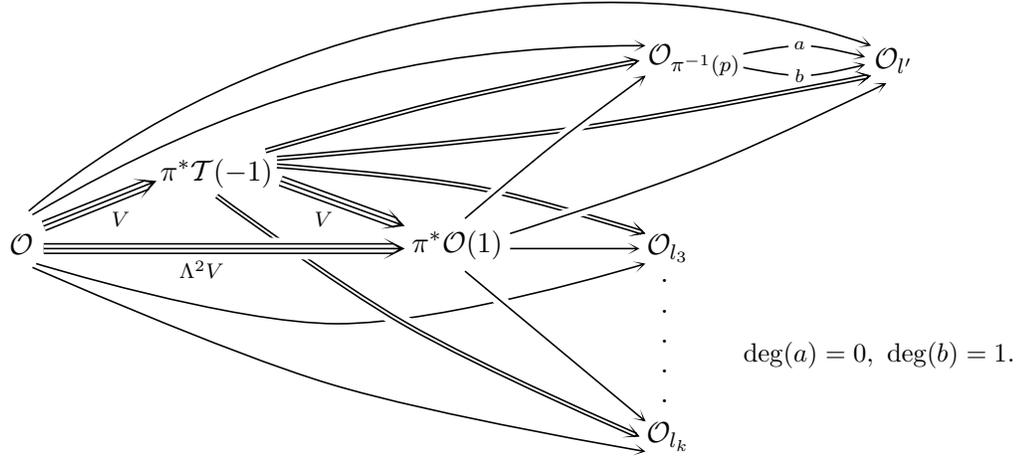
\begin{figure}[t]
\setlength{\unitlength}{1mm}
\begin{picture}(115,62)(30,-27)
\psset{unit=\unitlength,arrowsize=4pt 1,arrowlength=1,arrowinset=0.75,linewidth=0.6pt}
\put(29,0.5){\makebox(0,0)[cc]{$\O$}}
\put(55,10){\makebox(0,0)[cc]{$\pi^*\T(-1)$}}
\put(87,0.5){\makebox(0,0)[cc]{$\pi^*\O(1)$}}
\put(118.5,25){\makebox(0,0)[cc]{$\O_{\pi^{-1}(p)}$}}
\put(145,25){\makebox(0,0)[cc]{$\O_{l'}$}}
\put(115,0){\makebox(0,0)[cc]{$\O_{l_3}$}}
\multiput(114,-4.5)(0,-4){5}{.}
\put(115,-25){\makebox(0,0)[cc]{$\O_{l_k}$}}
\psbezier{->}(30.5,4.5)(57,20)(80,27)(112,27)
\psbezier{->}(30,5)(60,30)(105,40)(142,27)
\pscurve{->}(31,-2)(71,-10)(112,-2)
\pscurve[doubleline=true,border=2pt]{->}(55,7)(80,-10)(111.2,-25)
\psline[doubleline=true,doublesep=5\pslinewidth,border=2pt]{-}(32,0)(78.5,0)
\psline[linewidth=\pslinewidth,arrowsize=4pt 7]{->}(32,0)(80,0)
\psline[doubleline=true,doublesep=5\pslinewidth]{-}(32,3)(45.75,8.5)
\psline[linewidth=\pslinewidth,arrowsize=4pt 7]{->}(32,3)(47,9)
\psline[doubleline=true,doublesep=5\pslinewidth]{-}(63.5,9)(78.8,3.436)
\psline[linewidth=\pslinewidth,arrowsize=4pt 7]{->}(63.5,9)(80,3)
\pscurve[doubleline=true]{->}(61.5,13)(85,19.5)(111.2,25)
\pscurve[doubleline=true]{->}(63,12)(103,16)(142,23)
\pscurve[doubleline=true]{->}(63,11)(90,8)(112,2)
\pscurve[border=2pt]{->}(88,4)(100,14)(112,23)
\pscurve[border=2pt]{->}(94,2)(120,11)(144,22)
\psline{->}(94,0)(111.2,0)
\pscurve[border=2pt]{->}(88,-3)(100,-13)(112,-23)
\pscurve{->}(30.5,-2.5)(70,-18)(112,-27)
\put(50,-4){$\scriptstyle \Lambda^2 V$}
\put(41,3){$\scriptstyle V$}
\put(68,3){$\scriptstyle V$}
\pscurve{->}(125,25.9)(133,27)(141.2,25.5)
\pscurve{->}(125,24.1)(133,23)(141.2,24.5)
\psframe[fillstyle=solid,fillcolor=white,linecolor=white](131,27.5)(134,22.5)
\put(132.5,27){\makebox(0,0)[cc]{\tiny $a$}}
\put(132.5,23){\makebox(0,0)[cc]{\tiny $b$}}
\put(125,-15){\small $\deg(a)=0,\ \deg(b)=1.$}
\end{picture}
\caption{The cohomology algebra of the DG-quiver $\Bc_K(\tau)$ for a blowup
of $\PP^2$ with two infinitely close points.}
\label{quiver'}
\end{figure}

\begin{prop}
Let $X_K$ be the blowup of $\PP^2$ at a subscheme $K$ supported at a set of
$k-1$ points $\{ p, p_3,\dots, p_k\}$ and with length $2$ at the point $p$.
Then the sequence
\begin{equation}\label{excol:deg}
\tau=\left(\O_{X_K},\pi^* \T_{\PP^2}(-1), \pi^*\O_{\PP^2}(1), \O_{\pi^{-1}(p)}, \O_{l'}, \O_{l_3},\dots,\O_{l_k}\right)
\end{equation}
is a full exceptional collection on $X_K.$
\end{prop}

As before we can see that the sheaves $\O_{l_i}$ and $\O_{l_j}$ are
mutually orthogonal for all $i\ne j,$ and each $\O_{l_i}$ is orthogonal to both
$\O_{l'}$ and $\O_{\pi^{-1}(p)}.$
However, the collection $\tau$ is  not strong, because there are non-trivial morphisms from $\O_{\pi^{-1}(p)}$ to $\O_{l'}$
in degrees 0 and 1.
More precisely,
$$
\Hom(\O_{\pi^{-1}(p)}, \O_{l'})\cong \kk
\qquad
\text{and}
\qquad
\Ext^{1}(\O_{\pi^{-1}(p)}, \O_{l'})\cong \kk.
$$

Denote by $a$ and $b$ two morphisms from $\O_{\pi^{-1}(p)}$ to $\O_{l'}$ of
degrees 0 and 1 respectively.
It is easy to see that composition with the morphism $a$ gives isomorphisms between
the spaces $\Hom(F, \O_{\pi^{-1}(p)})$ and $\Hom(F, \O_{l'})$ for any
element $F$ of the exceptional collection $\tau$ (see Figure \ref{quiver'}).

Two approaches can be used to obtain an analogue of equivalence (\ref{equiv:del})
for this situation. The first possibility is to associate to the non-strong
exceptional collection $\tau$ a
differential graded algebra of homomorphisms,
and obtain an equivalence between the derived category of coherent sheaves
and the derived category of finitely generated (right) DG-modules over the DG-algebra of homomorphisms
of the exceptional collection. (One could also try to work in the framework of
$A_\infty$-algebras, which might be more appropriate here considering that
the mirror situation involves an $A_\infty$-category with non-zero $m_3$,
see \S \ref{ss:simpdegsympl}).

Another approach is to change the exceptional collection $\tau$ to another one
which is strong. There are natural operations on exceptional collections which
are called mutations and which allow us to obtain new exceptional collections
starting from a given one.

We omit the definition of mutations, which is classical and can be found in
many places. However, we note that the
left mutation of the exceptional collection (\ref{excol:deg}) in the pair
$(\pi^*\O_{\PP^2}(1), \O_{\pi^{-1}(p)})$ gives us a new exceptional collection
\begin{equation}\label{excolnew:deg}
\tau'=\left( \O_{X_K}, \pi^* \T_{\PP^2}(-1), \mathcal{M}, \pi^*\O_{\PP^2}(1),  \O_{l'}, \O_{l_3},\dots,\O_{l_k}\right)
\end{equation}
where $\M$ is the line bundle on $X_K$ which is the kernel of the surjection
$\pi^* \O_{\PP^2}(1)\to \O_{\pi^{-1}(p)}.$ This new exceptional collection
$\tau'$ is strong.

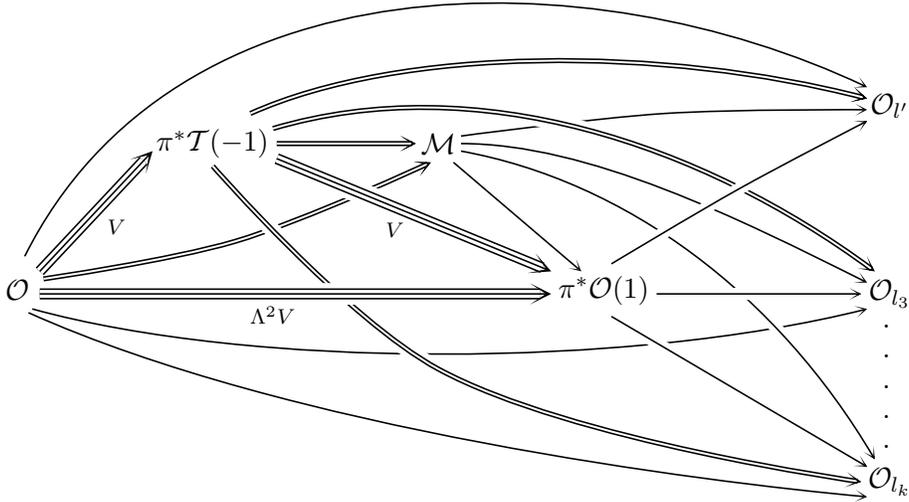
\begin{figure}[b]
\setlength{\unitlength}{1mm}
\begin{picture}(115,65)(30,-27)
\psset{unit=\unitlength,arrowsize=4pt 1,arrowlength=1,arrowinset=0.75,linewidth=0.6pt}
\put(29,0.5){\makebox(0,0)[cc]{$\O$}}
\put(55,20){\makebox(0,0)[cc]{$\pi^*\T(-1)$}}
\put(85,20){\makebox(0,0)[cc]{$\M$}}
\put(107,0.5){\makebox(0,0)[cc]{$\pi^*\O(1)$}}
\put(145,25){\makebox(0,0)[cc]{$\O_{l'}$}}
\put(145,0){\makebox(0,0)[cc]{$\O_{l_3}$}}
\multiput(144,-4.5)(0,-4){5}{.}
\put(145,-25){\makebox(0,0)[cc]{$\O_{l_k}$}}
\psbezier{->}(30,5)(50,45)(105,45)(142,27.5)
\psbezier{->}(31,-2)(60,-10)(110,-10)(142,-2)
\psline[doubleline=true]{->}(63.5,20)(82,20)
\pscurve[doubleline=true]{->}(32.5,2)(58,7)(84,17.5)
\psbezier[border=2pt]{->}(88,19)(110,15)(130,0)(143,-22)
\pscurve[doubleline=true,border=2pt]{->}(55,17)(85,-10)(141.2,-25)
\psline[doubleline=true,doublesep=5\pslinewidth,border=2pt]{-}(32,0)(98.5,0)
\psline[linewidth=\pslinewidth,arrowsize=4pt 7]{->}(32,0)(100,0)
\psline[doubleline=true,doublesep=5\pslinewidth]{-}(32,3)(45.8,17.72)
\psline[linewidth=\pslinewidth,arrowsize=4pt 7]{->}(32,3)(47,19)
\psline[doubleline=true,doublesep=5\pslinewidth,border=2pt]{-}(63.5,18)(98.8,3.48)
\psline[linewidth=\pslinewidth,arrowsize=4pt 7]{->}(63.5,18)(100,3)
\pscurve{->}(88,21)(115,24)(142,24.5)
\psbezier{->}(88,20)(105,20)(125,10)(142,1.5)
\psline{->}(87,17.5)(104,3)
\psbezier[doubleline=true,border=2pt]{->}(63,22)(90,30)(120,20)(143,3)
\psbezier[doubleline=true]{->}(60,24)(80,33)(115,33)(142,26)
\pscurve[border=2pt]{->}(108,4)(125,14)(142,23)
\psline[border=2pt]{->}(114,0)(141.2,0)
\pscurve[border=2pt]{->}(108,-3)(125,-13)(142,-23)
\psbezier{->}(30.5,-2.5)(65,-18)(120,-25)(142,-27)
\put(60,-4){$\scriptstyle \Lambda^2 V$}
\put(41,8){$\scriptstyle V$}
\put(78,7.5){$\scriptstyle V$}
\end{picture}
\caption{The quiver $\Bc_K(\tau')$ for a blowup of $\PP^2$ with
two infinitely close points.}
\label{quiver''}
\end{figure}

In fact, we can also consider the same left mutation when the
blown up points are all distinct, and obtain in that case as well a strong
exceptional collection $$\sigma'=(\O_{X_K}, \pi^* \T_{\PP^2}(-1), 
\mathcal{M}, \pi^*\O_{\PP^2}(1),  \O_{l_2}, \O_{l_3},\dots,\O_{l_k}),$$
which behaves very much like $\tau'$. The
distinguishing feature of the case where we blow up the point $p$ twice
is that in the exceptional collection $\tau'$ the composition map
\begin{equation}\label{eq:comptrivial}
\Hom(\M,\pi^*\O_{\PP^2}(1))\otimes \Hom(\pi^*\O_{\PP^2}(1),\O_{l'})\lto
\Hom(\M,\O_{l'})
\end{equation}
is identically zero, whereas for $\sigma'$ (i.e., when the points of $K$ are
distinct) the
corresponding composition is non-trivial. In this sense, the mutation
allows us to give a simple description of the behaviour of the
category under the degeneration process where two points of $K$ converge
to each other. Namely, the algebra $B_K(\tau')$ of homomorphisms of the
exceptional collection $\tau'$ is obtained as a degeneration of the algebra
of homomorphisms of the exceptional collection $\sigma'$ in which
the composition (\ref{eq:comptrivial}) becomes zero.

\begin{prop}
Let $X_K$ be the blowup of \,$\PP^2$ at a subscheme $K$ supported at a set of
$k-1$ points $\{ p, p_3,\dots, p_k\}$ and with length $2$ at the point $p$.
Then there is an equivalence
$
\db{\coh(X_K)}\cong\db{\mod B_K(\tau')},
$
where $B_K(\tau')$ is the algebra of
homomorphisms of the exceptional collection $\tau'.$
\end{prop}

In this context, the surface $X_K$ can again be recovered from the category
$\Bc_{K}(\tau')$. Namely, the
points $p, p_3,\dots, p_k$ can be determined by the same method as above.
To recover $X_K$, we also have to determine the position of the point $p'$
on the exceptional
curve of the blowup of the point $p.$ This is equivalent to finding
a tangent direction at the point $p.$
Consider  the kernel of the composition map
$$
\Hom(\O,\M)\otimes \Hom(\M, \O_{l'})\lto \Hom(\O, \O_{l'}).
$$
It is a one-dimensional subspace of $\Hom(\O, \M).$ The image of this subspace
in the space $V^*=\Hom(\O, \pi^*\O(1))$ determines a line in the projective
space $\PP(V)$ which passes through the point $p$ and hence a tangent
direction at $p$.

\subsection{Noncommutative deformations of Del Pezzo surfaces}
\label{ss:noncomm}

As before, let $X_K$ be the blowup of the projective plane at a set
$K$ of $k$ distinct points. Consider the strong exceptional collection
$$
\sigma=(\O, \pi^*\T_{\PP^2}(-1), \pi^*\O_{\PP^2}(1), \O_{l_1},\dots, \O_{l_k}).
$$
By the discussion in \S \ref{ss:delpezzo}, the derived category
of coherent sheaves $\db{\coh(X_K)}$ is equivalent to the category of
finitely generated (right) modules
over the algebra $\Be_K$ of homomorphisms of $\sigma$.
The algebra $\Be_K$ can also be represented by the category $\Bc_K$ associated to
the exceptional collection $\sigma$ (see Figure \ref{quiver}).

The category $\Bc_K$ has strictly more deformations than the surface $X_K.$
We saw above that the surface $X_K$ can be reconstructed from the category
$\Bc_K$, and that
the deformation of the surface $X_K$ is controlled by the variation of the set
$K\subset\PP^2.$

A general deformation of the category $\Bc_K$ can be viewed as
the category of an exceptional collection on a noncommutative deformation
of the surface $X_K.$ In other words, if $\Bc_{K,\mu}$ is a deformation of the
category $\Bc_K$ then the bounded derived category $\db{\mod\Be_{K,\mu}}$ of
finitely generated (right) modules over the algebra $\Be_{K,\mu}$ will be
viewed as the derived category
of coherent sheaves on a noncommutative surface $X_{K,\mu}.$
Any such noncommutative surface can be represented as the blowup of a noncommutative
plane $\PP^2_{\mu}$ at some set $K$ consisting of $k$ of its ``points''.
This procedure is discussed in detail in \cite{VdB}.

In the rest of this section, we describe the deformations of the category $\Bc_K.$
Recall that a deformation of a category is, by definition, a deformation of
its composition law. We proceed in two steps.
The first step is to describe the deformations of the subcategory $\Bc(\sigma_0)$
associated to the subcollection $\sigma_0=(\O, \pi^*\T(-1), \pi^*\O(1)).$
This subcategory $\Bc(\sigma_0)$ is the category of homomorphisms
of the full strong exceptional collection $(\O, \T(-1), \O(1))$ on
$\PP^2.$  Therefore, considering a deformation of the subcategory
$\Bc(\sigma_0)$ we obtain a noncommutative deformation $\PP^2_{\mu}$ of the
projective plane.
The second step is to describe the deformations of all other compositions in the category $\Bc_K.$
These deformations correspond to variations of the set of ``points'' $K$ on the
noncommutative projective plane $\PP^2_{\mu}.$

\begin{figure}[t]
\xy
\POS (-30,0) *+{}="0"
\POS (30,0)*+{F_0}="e0"
\POS (50,15)*+{F_1}="e1"
\POS (70,0)*+{F_2}="e2"

\POS"e0" \ar@3^{U}@<0.5ex> "e1"
\POS"e0" \ar@3_{W} "e2"
\POS"e1" \ar@3^{V}@<0.2ex> "e2"

\endxy
\caption{The quiver $\Bc_{\mu}$ for a noncommutative $\PP^2_{\mu}.$}
\label{quiver:small}
\end{figure}

Noncommutative deformations of the projective plane have been described in \cite{ATV, BP}.
Any deformation of the category $\Bc(\sigma_0)$ is a category with three ordered objects
$F_0, F_1, F_2$
and with three-dimensional spaces of homomorphisms from  $F_i$ to $F_j$ when $i<j$
(see Figure \ref{quiver:small}). Any such category $\Bc_{\mu}$ is determined by
the composition tensor $\mu: V\otimes U\to W.$ We will consider only the nondegenerate (geometric)
case, where the restrictions $\mu_u=\mu(\cdot,u): V\to W$ and 
$\mu_v=\mu(v,\cdot): U\to W$ have rank at least two for all nonzero elements
$u\in U$ and $v\in V$, and the composition of $\mu$ with any nonzero linear
form on $W$ is a bilinear form of rank at least two on $V\otimes U$.
The equations $\det \mu_u=0$ and $\det\mu_v=0$ define closed subschemes
$\Gamma_U\subset\PP(U)$ and $\Gamma_V\subset\PP(V).$
Namely, up to projectivization the set of points of $\Gamma_U$ 
(resp.\ $\Gamma_V$) consists of all
$u\in U$ (resp.\ $v\in V$) for which the rank of $\mu_u$ (resp.\ $\mu_v$) is equal to $2.$
It is easy to see that the correspondence which associates to
a vector $v\in V$ the kernel of the map $\mu_v: U\to W$ defines
an isomorphism between $\Gamma_V$ and $\Gamma_U$.
Moreover, under these circumstances $\Gamma_V$ is either the entire
projective plane $\PP(V)$ or a cubic
in $\PP(V).$ If $\Gamma_V=\PP(V)$ then $\mu$ is isomorphic to the tensor $V\otimes V\to \Lambda^2 V,$
i.e. we get the usual projective plane $\PP^2.$

Thus, the non-trivial case is the situation where $\Gamma_V$ is a cubic, i.e.\
an elliptic curve which we now denote by $E$.
This elliptic curve comes equipped with two embeddings into the
projective planes $\PP(U)$ and $\PP(V)$ respectively; by restriction of
$\O(1)$ these embeddings
determine two line bundles $\L_1$ and $\L_2$ of degree $3$ over $E$,
and it can be checked that $\L_1\ne \L_2.$ This construction has a converse:

\begin{construction}\label{constr:mu}
The tensor $\mu$ can be reconstructed from the triple
$(E, \L_1, \L_2).$
Namely, the spaces $U, V$ are isomorphic to
$H^0(E, \L_1)^*$ and $H^0(E, \L_2)^*$ respectively, and the space
$W$ is dual to the kernel of the canonical map
$$
H^0(E, \L_1) \otimes H^0(E, \L_2)\lto H^0(E, \L_1\otimes\L_2),
$$
which induces the tensor $\mu: V\otimes U\lto W$.
\end{construction}

The details of these constructions and statements can be found in \cite{ATV,BP}.

\begin{remark}\label{rmk:comm}
Note that  we can also consider  a triple $(E, \L_1, \L_2)$ such that
$\L_1\cong \L_2.$ Then the procedure described above produces
a tensor $\mu$ with $\Gamma_V\cong \PP(V)$, which defines the usual
commutative projective plane. Thus, in this particular case the tensor  $\mu$ does not depend on the curve $E.$
\end{remark}

Now we describe the deformations of the other compositions in the category $\Bc_K$.
Given a category $\Bc_{\mu}$ of the form described above, corresponding to a
noncommutative projective plane $\PP^2_{\mu}$, and given a
set $K=\{ p_1,\dots,p_k \}$ of $k$ points on the elliptic curve $E$, 
we can construct a category $\Bc_{K,\mu}$ in the following manner.
A point $p_j\in E\subset\PP(U)$ determines a one-dimensional subspace of
$U$, generated by a vector $u_j\in U$. The map $\mu_{u_j}: V\to W$
has rank $2$; denote by $v_j$ a non-zero vector in its kernel.
The category $\Bc_{K,\mu}$ is then constructed
from the category $\Bc_{\mu}$ by adding $k$ mutually orthogonal
objects $\O_{l_j}$ for $j=1,\dots, k$,
and defining the spaces of morphisms by the rule
$$
\Hom(F_2, \O_{l_j})=\kk,
\qquad
\Hom(F_1, \O_{l_j})= V/\Ker \mu_{u_j}=V/\langle v_j\rangle,
\qquad
\Hom(F_0, \O_{l_j})= W/\Im \mu_{v_j}.
$$
The two composition tensors involving $\Hom(F_2,\O_{l_j})$ are defined in
the obvious manner as suggested by the notation. The only non-obvious
composition is the map $V/\langle v_j\rangle \otimes U\to W/\Im\mu_{v_j}$, which
is by definition induced by the tensor $\mu: V\otimes U\lto W.$

Conversely, if we consider a category $\Bc_{K,\mu}$ which is a
deformation of $\Bc_K$ and an extension of the category
$\Bc_{\mu}$, then the kernel of  the composition map
$$
\Hom(F_2, \O_{l_j})\otimes V\lto \Hom(F_1, \O_{l_j})
$$
defines a one-dimensional subspace $\langle v_j\rangle \subset V$.
The map $\mu_{v_j}$ must have rank $2$, since otherwise $\mu_{v_j}$
would be an isomorphism and the composition map $\Hom(F_2,
\O_{l_j})\otimes W\lto \Hom(F_0, \O_{l_j})$ would vanish identically, which
by assumption is not the case.
Therefore, the objects $\O_{l_j}$  correspond to points on the curve $E$.

Thus, any category $\Bc_{K,\mu}$ is defined by the data
$(E, \L_1, \L_2, p_1,\dots, p_k),$ where
$E$ is a cubic, $\L_1, \L_2$ are line bundles of degree $3$ on $E$, and
$p_1,\dots, p_k$ is a set of distinct points on $E$.
If $\L_1\cong\L_2,$ then we obtain the category $\Bc_{K}$ related to a 
blowup of the usual commutative projective plane. In the general case,
the bounded derived category $\db{\mod\Be_{K,\mu}}$ of finite rank
modules over the algebra $\Be_{K,\mu}$ is viewed as the derived
category of coherent sheaves
on the non-commutative surface $X_{K,\mu}$, which is a blowup of
$k$ points on the non-commutative projective plane $\PP^2_{\mu}$.

\medskip

A standard approach  to noncommutative geometry is to determine a
noncommutative variety either by an abelian category of
(quasi)coherent sheaves on it or by a noncommutative (graded)
algebra which is considered as its (homogeneous) coordinate ring. The question
of how to define the abelian category of coherent sheaves on Del
Pezzo surfaces and on other blowups of surfaces is discussed in the
paper \cite{VdB}.
We briefly describe one of the possible approaches. It is very
important to note that the category $\db{\mod\Be_{K,\mu}}$
possess a Serre functor $S$, i.e.\ an additive autoequivalence
for which  there are bi-functorial isomorphisms
$$
\Hom(X, SY)\stackrel{\sim}{\lto} \Hom(Y, X)^*
$$
for any $X, Y\in \db{\mod\Be_{K,\mu}}$. In the case
of the bounded derived  category of finite rank modules
over a finite dimensional algebra of finite homological dimension,
the Serre functor is the  functor which takes a
complex of modules $M^{\bdot}$ to the complex
$\bR\Hom_{\Be_{K,\mu}}(M, \Be_{K,\mu})^*$. The Serre functor
is an exact autoequivalence.

Now we can take the projective module $P_0$ (corresponding to $\O$, see
the discussion after Definition \ref{def:algexc}) and consider the
sequence of objects $R_m=S^{m}[-2m] P_0$ for all $m\in \Z.$ Let
us consider the subcategory $\A\subset\db{\mod\Be_{K,\mu}}$
consisting of all objects $F$ such that
$$
\Hom(R_m, F[i])=0
\quad
\text{for all}
\ i\ne 0
\ \text{and sufficiently large}
\ m\gg 0.
$$
If the category $\A$ is abelian and its bounded derived
category $\db{\A}$ is equivalent to $\db{\mod\Be_{K,\mu}}$
then $\A$ can be considered as the category of coherent
sheaves on the noncommutative surface $X_{K,\mu}$, and $X_{K,\mu}$
can be called a noncommutative Del Pezzo surface.

The reason of such a definition of the abelian category of coherent sheaves
on a noncommutative Del Pezzo surface is inspired by the commutative case.
In the commutative case the Serre functor is isomorphic to the functor
$\otimes\O(K)[2],$ where $\O(K)$ is the canonical line bundle.
Hence, for usual commutative surfaces the objects $R_m$ are isomorphic to the
invertible sheaves $\O(mK).$ Since for a Del Pezzo surface $X$ the anticanonical sheaf $\O(-K)$ is ample,
we have $H^i(X, F(-mK))=0$ for all $i\ne 0$ and any coherent sheaf $F$ when $m$ is sufficiently large.
This property makes it possible to separate pure coherent sheaves from other complexes of coherent sheaves.

We can also consider
the graded space $A=\bigoplus_{p=0}^{\infty} \Hom(R_0, R_{-p})$ and  can
endow it with the structure of a graded algebra using the isomorphisms
$\Hom(R_0, R_{-p})\cong\Hom(R_{i}, R_{i-p})$ given by the
functors $S^{i}[-2i]$ for all $i\in \Z.$
This algebra can be considered as the homogeneous coordinate ring
of a noncommutative Del Pezzo surface. It seems that such rings  are  noncommutative deformations of
homogeneous commutative coordinate rings of usual Del Pezzo surfaces.

In any case,
these remarks about abelian categories of coherent sheaves on noncommutative
Del Pezzo surfaces will not be needed in the rest of the argument. We will only
use the description of the bounded derived category of coherent sheaves on
the noncommutative blowup
$X_{K,\mu}$ in terms of finite rank modules over the algebra $\Be_{K,\mu},$
i.e.\ we state an equivalence of triangulated categories
\begin{equation}
\db{\coh(X_{K,\mu})}\cong \db{\mod\Be_{K,\mu}}.
\end{equation}


\section{The mirror Landau-Ginzburg models}\label{s:mirrors}

\subsection{Compactification of the mirror of $\CP^2$}\label{ss:compactify}

As mentioned in the introduction, the mirror of $\CP^2$ is 
an elliptic fibration with 3 singular
fibers, determined by (a fiberwise compactification of) the superpotential
$W_0=x+y+{1}/{xy}$ on $(\C^*)^2$. This Landau-Ginzburg model
compactifies naturally to an elliptic fibration
$\ol{W_0}:\ol{M}\to\CP^1$, which we now describe.

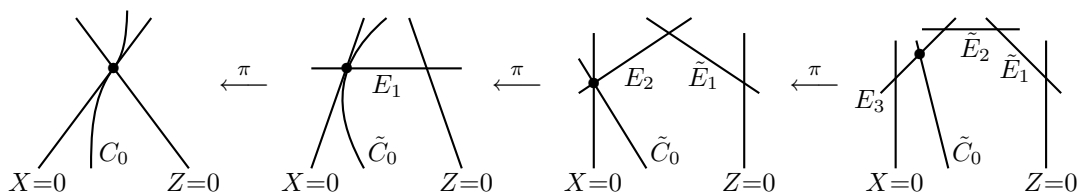
\begin{figure}[t]
\setlength{\unitlength}{1cm}
\begin{picture}(3.5,2.3)(-1.5,-1.2)
\psset{unit=\unitlength}
\psline(-1,-1)(0.5,1)
\psline(1,-1)(-0.5,1)
\put(0,0.333){\circle*{0.15}}
\psbezier(-0.3,-1)(-0.3,-0.6)(-0.3,-0.066)(0,0.333)
\psbezier(0,0.333)(0.125,0.5)(0.17,0.7)(0.18,1.1)
\put(-1.4,-1.3){\small $X$=0}
\put(0.7,-1.3){\small $Z$=0}
\put(-0.2,-0.9){\small $C_0$}
\put(1.4,0){\large $\stackrel{\pi}{\longleftarrow}$}
\end{picture}
\begin{picture}(3.5,2.3)(-1.5,-1.2)
\psset{unit=\unitlength}
\psline(-1,-1)(-0.3,1)
\psline(1,-1)(0.3,1)
\psline(-1,0.33)(1,0.33)
\put(-0.53,0.33){\circle*{0.15}}
\psbezier(-0.3,-1)(-0.63,-0.4)(-0.83,0.3)(0,1)
\put(-1.4,-1.3){\small $X$=0}
\put(0.7,-1.3){\small $Z$=0}
\put(-0.25,-0.9){\small $\tilde{C}_0$}
\put(-0.2,0){\small $E_1$}
\put(1.4,0){\large $\stackrel{\pi}{\longleftarrow}$}
\end{picture}\ \
\begin{picture}(3.5,2.3)(-1.5,-1.2)
\psset{unit=\unitlength}
\psline(-1,-1)(-1,0.8)
\psline(1,-1)(1,0.8)
\psline(-1.2,0)(0.3,1)
\psline(1.2,0)(-0.3,1)
\put(-1,0.133){\circle*{0.15}}
\psline(-0.3,-1)(-1.2,0.46)
\put(-1.4,-1.3){\small $X$=0}
\put(0.7,-1.3){\small $Z$=0}
\put(-0.25,-0.9){\small $\tilde{C}_0$}
\put(0.25,0.1){\small $\tilde{E}_1$}
\put(-0.6,0.1){\small $E_2$}
\put(1.6,0){\large $\stackrel{\pi}{\longleftarrow}$}
\end{picture}
\quad
\begin{picture}(3,2.3)(-1.5,-1.2)
\psset{unit=\unitlength}
\psline(-1,-1)(-1,0.7)
\psline(1,-1)(1,0.7)
\psline(-1.2,0)(-0.2,1)
\psline(1.2,0)(0.2,1)
\psline(-0.65,0.85)(0.65,0.85)
\put(-0.68,0.52){\circle*{0.15}}
\psline(-0.3,-1)(-0.725,0.7)
\put(-1.4,-1.3){\small $X$=0}
\put(0.7,-1.3){\small $Z$=0}
\put(-0.25,-0.9){\small $\tilde{C}_0$}
\put(0.38,0.20){\small $\tilde{E}_1$}
\put(-0.15,0.48){\small $\tilde{E}_2$}
\put(-1.55,-0.2){\small $E_3$}
\end{picture}
\caption{The successive blowups at $(0:1:0)$.}
\label{fig:3blowups}
\end{figure}

Compactifying $(\C^*)^2$ to $\CP^2$, we can view $W_0$ as the
quotient of the two homogeneous degree 3 polynomials $P_0=X^2Y+XY^2+Z^3$ and
$P_\infty=XYZ$, which define a pencil of cubics with three base points of
multiplicities respectively $4$, $4$, and $1$. Namely, the cubic $C_0$
defined by $P_0$ intersects the line $X=0$ at $(0:1:0)$ (with multiplicity
3), the line $Y=0$ at $(1:0:0)$ (with multiplicity 3), and the line $Z=0$ at
$(0:1:0)$, $(1:0:0)$ and $(1:-1:0)$. Blow up $\CP^2$ three times
successively at the point where the cubic $C_0$ and the line $X=0$
(or their proper transforms) intersect each other, i.e.\ first at the point
$(0:1:0)$, and then twice at suitable points of the exceptional divisors
(see Figure \ref{fig:3blowups}).
Similarly, blow up three times the intersection of the cubic $C_0$ with
the line $Y=0$.

 Let $\tilde{C}_0$ be the proper transform of $C_0$ under
these blowups, and let $\tilde{C}_\infty$ be the configuration of 9 rational
curves formed by the proper transforms of the three coordinate lines and the
exceptional divisors of the six blowups
(so, in Figure \ref{fig:3blowups}, all components other than
$\tilde{C}_0$ are eventually part of $\tilde{C}_\infty$).
Then $\tilde{C}_0$ and
$\tilde{C}_\infty$ intersect transversely at three smooth points, and define a
pencil of elliptic curves representing the anticanonical class in $\CP^2$
blown up six times. The complement of $\tilde{C}_\infty$ identifies
with $(\C^*)^2$, and the restriction of the $\CP^1$-valued map defined by
the pencil to this open subset coincides with $W_0$.
Blowing up the three points where $\tilde{C}_0$ and $\tilde{C}_\infty$ intersect,
we obtain a rational elliptic surface $\ol{M}$, and the pencil becomes
an elliptic fibration $\ol{W_0}:\ol{M}\to\CP^1$, which provides a natural
compactification of $W_0:(\C^*)^2\to\C$.

\begin{figure}[b]
\setlength{\unitlength}{1cm}
\begin{picture}(5,3)(-2.1,-0.8)
\psset{unit=\unitlength}
\put(-2.1,-0.5){\line(1,0){5}}
\put(-1.5,-0.5){\circle*{0.1}}
\put(-0.75,-0.5){\circle*{0.1}}
\put(0,-0.5){\circle*{0.1}}
\put(2,-0.5){\circle*{0.1}}
\put(-1.75,-0.95){$3j^2$}
\put(-1.5,1){\circle*{0.1}}
\put(-0.95,-0.95){$3j$}
\put(-0.75,1){\circle*{0.1}}
\put(-0.1,-0.95){$3$}
\put(0,1){\circle*{0.1}}
\put(1.85,-0.95){$\infty$}
\put(1.65,0.12){\circle*{0.1}}
\put(1.3,0.5){\circle*{0.1}}
\put(1.15,1.1){\circle*{0.1}}
\put(1.53,1.86){\circle*{0.1}}
\put(2.35,0.12){\circle*{0.1}}
\put(2.7,0.5){\circle*{0.1}}
\put(2.85,1.1){\circle*{0.1}}
\put(2.47,1.86){\circle*{0.1}}
\put(2,2.03){\circle*{0.1}}
\pscurve(-1.5,0)(-1.5,1)(-1.8,1.3)(-1.8,0.7)(-1.5,1)(-1.5,2)
\pscurve(-0.75,0)(-0.75,1)(-1.05,1.3)(-1.05,0.7)(-0.75,1)(-0.75,2)
\pscurve(0,0)(0,1)(-0.3,1.3)(-0.3,0.7)(0,1)(0,2)
\psline{-}(1.5,0.1)(2.5,0.1)
\psline{-}(1.8,-0.05)(1.2,0.6)
\psline{-}(1.35,0.3)(1.1,1.3)
\psline{-}(1.1,1)(1.6,2)
\psline{-}(1.4,1.8)(2.2,2.1)
\psline{-}(2.2,-0.05)(2.8,0.6)
\psline{-}(2.65,0.3)(2.9,1.3)
\psline{-}(2.9,1)(2.4,2)
\psline{-}(2.6,1.8)(1.8,2.1)
\end{picture}
\caption{The singular fibers of $\ol{W_0}$.}\label{fig:W0crit}
\end{figure}
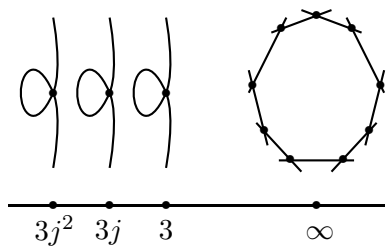

The meromorphic function $\ol{W_0}$ has 12 isolated non-degenerate critical
points. Three of them are the pre-images of the points $(1:1:1)$, $(j:j:1)$,
and $(j^2:j^2:1)$ ($j=e^{2i\pi/3}$), and correspond to the three
critical points of $W_0$ in $(\C^*)^2$ (with associated critical values
$3$, $3j$, and $3j^2$). The nine other critical points all lie in the fiber
above infinity: they are the nodes of the reducible configuration
$\tilde{C}_\infty$ (see Figure \ref{fig:W0crit}).
\medskip

This compactification process can also be described in a more symmetric
manner by viewing $(\C^*)^2$ as the surface $\{xyz=1\}\subset(\C^*)^3$, and
$W_0=x+y+z$. Compactifying $(\C^*)^3$ to $\CP^3$ leads one to consider
the cubic surface
$\{XYZ=T^3\}\subset\CP^3$, which presents $A_2$ singularities at the
three points $(1:0:0:0)$, $(0:1:0:0)$, and $(0:0:1:0)$. After blowing up $\CP^3$
at these three points, we obtain a smooth cubic surface, in which the
hyperplane sections $\tilde{C}_0=\{X+Y+Z=0\}$  and
$\tilde{C}_\infty=\{T=0\}$ define a pencil of elliptic curves
with three base points.
As before, $\tilde{C}_0$ is a smooth elliptic curve, and $\tilde{C}_\infty$ is
a configuration of 9 rational curves (the proper transforms of the three
coordinate lines where the singular cubic surface intersects the plane
$T=0$, and the six $-2$-curves arising from the resolution of the
singularities). Blowing up the three points of $\tilde{C}_0\cap
\tilde{C}_\infty$, we again obtain a rational elliptic surface, and an elliptic
fibration with 12 isolated critical points, 9 of which lie in the fiber
above infinity (as in Figure \ref{fig:W0crit}).

\subsection{The vanishing cycles of $\ol{W_0}$}\label{ss:critw0}

Each singular fiber of $\ol{W_0}$ is obtained from the regular fiber by
collapsing a certain number of {\it vanishing cycles}, and the monodromy of
the fibration around a singular fiber is given by a product of Dehn twists
along these vanishing cycles. In this section, we determine the homology
classes of the vanishing cycles associated to the critical points of
$\ol{W_0}$.

More precisely, consider the fiber $\Sigma_0=\ol{W_0}{}^{-1}(0)$, which is
a smooth elliptic curve (in fact, the proper transform of the curve
called $\tilde{C}_0$ in \S \ref{ss:compactify}), and consider the following
ordered collection of arcs
$(\gamma_i)_{0\le i\le 3}$ joining the origin to the various critical values
of $\ol{W_0}$: $\gamma_0,\,\gamma_1,\,\gamma_2$ are
straight line segments joining the origin to $\lambda_0=3$, $\lambda_1=3j^2$,
and $\lambda_2=3j$ respectively, and $\gamma_3$ is the straight line
$e^{i\pi/3}\R_+$ joining the origin to $\lambda_3=\infty$.

Using parallel transport (with respect to an arbitrary horizontal
distribution) along the arc $\gamma_i$, we can associate a vanishing
cycle to each critical point $p\in \ol{W_0}{}^{-1}(\lambda_i)$; this
vanishing cycle is well-defined up to isotopy, and in particular we can
consider its homology class in $H_1(\Sigma_0,\Z)\simeq \Z^2$ (well-defined
up to a choice of orientation). If we fix
a symplectic structure on $\ol{M}$ for which the fibers of $\ol{W_0}$ are
symplectic submanifolds, then we have a canonical horizontal distribution
(given by the symplectic orthogonal to the fiber), which allows us to
consider the vanishing cycles as Lagrangian submanifolds of $\Sigma_0$,
well-defined up to {\it Hamiltonian isotopy}; in \S \ref{s:lagvc} this
will be of utmost importance, but for now we ignore the symplectic structure
and only view $\ol{W_0}$ as a topological fibration.

\begin{lemma}\label{l:vch1}
In terms of a suitable basis $\{a,b\}$ of $H_1(\Sigma_0,\Z)$, the
vanishing cycles $L_0,L_1,L_2$ associated to the critical values
$\lambda_0,\lambda_1,\lambda_2$ $($and the arcs $\gamma_0,\gamma_1,\gamma_2)$
represent the classes $[L_0]=-2a-b$, $[L_1]=-a+b$, and $[L_2]=a+2b$,
respectively;
and the vanishing cycles $L_3,\dots,L_{11}$ associated to the nine critical
points in the fiber at infinity represent the class $[L_3]=\dots=[L_{11}]=
a+b$.
\end{lemma}

\proof
The vanishing cycles $L_0,L_1,L_2$
are exactly those of the mirror of $\CP^2$, and are well-known (cf.\ e.g.\
\cite{Se2} or \cite{AKO}). In particular it is known that, choosing a
suitable homology basis $\{a,b\}$ for $H_1(\Sigma_0,\Z)$, and fixing
appropriate orientations of $L_0,L_1,L_2$, we have
$[L_0]=-2a-b$, $[L_1]=-a+b$, and $[L_2]=a+2b$ (cf.\ e.g.\ Figure 14 in
\cite{AKO}).

We now consider the 9 critical points in the fiber at infinity.
It is clear that
$L_3,\dots,L_{11}$ admit disjoint representatives, and hence are all
homologous. Their homology class can be determined by considering the
monodromy of the elliptic fibration $\ol{W_0}$, which is given by the
product of the positive Dehn twists along the vanishing cycles. Considering
the action on $H_1(\Sigma_0,\Z)$, and still using the basis $\{a,b\}$
considered above, the monodromies around the critical values
$\lambda_0,\lambda_1,\lambda_2$ are given by
$$\tau_0=\left(\begin{array}{rr}\!\!-1&4\\\!\!-1&3\end{array}\right),\quad
\tau_1=\left(\begin{array}{rr}2&1\\\!\!-1&0\end{array}\right),\quad\mathrm{and}
\ \tau_2=\left(\begin{array}{rr}\!-1&1\\\!-4&3\end{array}\right),
$$ while the monodromy around the fiber at infinity is given by $\tau^9$,
where $\tau$ is the positive Dehn twist along $[L_3]=\dots=[L_{11}]$.
On the other hand, because the arcs $\gamma_0,\dots,\gamma_3$ are ordered
clockwise around the origin, we have $\tau_0\tau_1\tau_2\tau^9=1$.
Therefore,
$$\tau^9=\left(\begin{array}{rr}\!-8&9\\\!-9&10\end{array}\right),$$
and considering $\mathrm{Ker}(\tau^9-1)$ we obtain
$[L_3]=\dots=[L_{11}]=a+b$.
\endproof

\proof[Alternative proof] (compare with \S 4.2 of \cite{AKO}).
Recall that $\ol{M}$ is obtained from
$\CP^2$ by successive blowups of the base points of the pencil of cubics
defined by $P_0=X^2Y+XY^2+Z^3$ and $P_\infty=XYZ$. Consider the ruled surface
$F$ obtained by blowing up $\CP^2$ just once at the point $(0\!:\!1\!:\!0)$:
the projection $(X\!:\!Y\!:\!Z)\mapsto (X\!:\!Z)$ naturally extends into a
fibration $\pi_x:F\to\CP^1$, of which the exceptional divisor is a section.
For $\lambda\in\CP^1$, denote by $\hat{C}_\lambda$ the proper transform of the
plane cubic $C_\lambda$ defined by $P_0-\lambda P_\infty$, which is also the image of
$\ol{W_0}{}^{-1}(\lambda)$ under the natural projection $p:\ol{M}\to F$.

The restriction $\pi_{x,\lambda}$ of $\pi_x$ to $\hat{C}_\lambda$ has degree two, and for
$\lambda\not\in\mathrm{crit}(\ol{W_0})$ its four branch points are
associated to distinct critical values in $\CP^1$, namely zero and the
three roots of the equation $x(\lambda-x)^2=4$. Indeed,
since $C_\lambda$ always has an order 3 tangency with the line $X=0$
at $(0\!:\!1\!:\!0)$, $\hat{C}_\lambda$ is always tangent to the fiber
$\pi_x^{-1}(0)$. The three other branch points are the critical points
of the projection to the first coordinate on $(\C^*)^2\cap C_\lambda=\{
(x,y)\in(\C^*)^2,\ xy(\lambda-x-y)=1\}$; viewing
$xy(\lambda-x-y)=1$ as a quadratic equation in the variable $y$,
the discriminant is $x(\lambda-x)^2-4$.

\begin{figure}[b]
\setlength{\unitlength}{1.6cm}
\begin{picture}(4,1.8)(-2,-0.9)
\psset{unit=\unitlength}
\put(1,0){\circle*{0.08}}
\put(-0.5,0.866){\circle*{0.08}}
\put(-0.5,-0.866){\circle*{0.08}}
\put(0,0){\circle*{0.08}}
\psbezier(-0.5,0.866)(0.5,0.288)(0.5,-0.288)(-0.5,-0.866)
\psbezier(-0.5,0.866)(-0.5,-0.288)(0,-0.576)(1,0)
\psbezier(-0.5,-0.866)(-0.5,0.288)(0,0.576)(1,0)
\psline(0,0)(-0.5,-0.866)
\psbezier(-0.5,0.866)(0.166,0.866)(0.666,0.555)(1,0)
\put(-0.66,-0.5){\makebox(0,0)[cc]{$\delta_2$}}
\put(-0.05,-0.8){\makebox(0,0)[cc]{$\delta_0$}}
\put(-0.66,0.5){\makebox(0,0)[cc]{$\delta_1$}}
\put(-0.32,-0.35){\makebox(0,0)[cc]{$\delta'$}}
\put(0.73,0.62){\makebox(0,0)[cc]{$\delta''$}}
\end{picture}
\caption{The projections of the vanishing cycles of $\ol{W_0}$}
\label{fig:projvcs}
\end{figure}
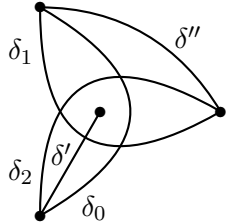

As $\lambda$ tends to $\lambda_i$ ($i\in\{0,1,2\}$), two of the
critical values of $\pi_{x,\lambda}$ converge to each other; keeping
track of the manner in which these critical values coalesce when $\lambda$
varies from $0$ to $\lambda_i$ along the arc $\gamma_i$, we obtain an arc
$\delta_i\subset\CP^1$, with end points in $\mathrm{crit}(\pi_{x,0})$
(see Figure~\ref{fig:projvcs}).
The lift of $\delta_i$ under the double cover $\pi_{x,0}$ is
(up to homotopy) the vanishing cycle $L_i$
(note that the projection $p:\ol{M}\to F$ allows us to implicitly
identify $\hat{C}_\lambda$ with
$\ol{W_0}{}^{-1}(\lambda)$ for $\lambda\neq\infty$).

Similarly, the behavior of the critical values of $\smash{\pi_{x,\lambda}}$ as
$\lambda$ tends to infinity describes the degeneration of $\hat{C}_\lambda$
to the singular
configuration $\hat{C}_\infty$, which consists of two sections and two
fibers of $\pi_x:F\to\CP^1$ (the fibers above $0$ and $\infty$, the
exceptional section, and the pre-image of the line $Y=0$). Namely, as
$\lambda$ tends to infinity along the arc $\gamma_3$, the critical value
with argument $-2\pi/3$ approaches zero, while the two other roots of
$x(\lambda-x)^2-4$
tend to infinity. The manner in which pairs of critical values coalesce is
encoded by the arcs $\delta'$ and $\delta''$ in Figure \ref{fig:projvcs},
and the four vanishing cycles associated to the degeneration are essentially
the lifts under $\pi_{x,0}$ of closed loops which bound regular
neighborhoods of the arcs $\delta'$ and $\delta''$; they all represent
the same homotopy class inside $\hat{C}_0$.

Recall that $\ol{W_0}{}^{-1}(\infty)\simeq \tilde{C}_\infty$ is obtained
from $\hat{C}_\infty$ by repeatedly blowing up two of the nodes.
Taking pre-images under these blowup operations,
the vanishing cycles associated to the two other nodes of $\hat{C}_\infty$
are naturally identified with two of the nine vanishing cycles
$L_3,\dots,L_{11}$ associated to the fiber at infinity of $\ol{W_0}$.
In particular, these vanishing cycles represent the same homology class
in $H_1(\Sigma_0,\Z)\simeq H_1(\hat{C}_0,\Z)$ as the lifts of
$\delta'$ and $\delta''$.

It is then easy to check that, for suitable choices of orientations, we have
$[L_0]\cdot[L_1]=[L_0]\cdot[L_2]=[L_1]\cdot[L_2]=-3$, $[L_0]\cdot[L_{3+i}]=
[L_2]\cdot[L_{3+i}]=-1$, and $[L_1]\cdot[L_{3+i}]=-2$, which completes
the proof of Lemma \ref{l:vch1}.
\endproof

\subsection{The vanishing cycles of $(M_k,W_k)$} \label{ss:critwk}

Recall from the introduction that our proposal for the mirror of a
Del Pezzo surface $X_K$
obtained from $\CP^2$ by blowing up $k\le 8$ generic points is an elliptic
fibration $W_k:M_k\to\C$, obtained by deforming the fibration $\ol{W_0}$ to
another elliptic fibration $\ol{W_k}:\ol{M}\to\CP^1$, and considering the
restriction to $M_k=\ol{M} \setminus \ol{W_k}{}^{-1}(\infty)$. More
precisely, remember that $\ol{W_k}$ has $3+k$ irreducible nodal fibers
corresponding to critical values $\lambda_0,\dots,\lambda_{k+2}\in\C$,
of which the first three correspond
naturally to the irreducible nodal fibers of $\ol{W_0}$, while the $k$
other finite critical values correspond to the deformation of critical
points in $\ol{W_0}{}^{-1}(\infty)$ towards finite values of the
superpotential. Meanwhile, $\ol{W_k}{}^{-1}(\infty)$ is a singular fiber
with $9-k$ components.

While the precise locations of the critical values $\lambda_i$ are closely
related to the complex structure on $M_k$, they are essentially irrelevant
from the point of view of symplectic topology and categories of vanishing
cycles. Indeed, if we consider a family
$(M_{k,t},W_{k,t})$ of fibrations indexed by a real parameter $t$, with
the property that for all $t$ the critical points of $W_{k,t}$ are
isolated and non-degenerate, then the vanishing cycles
remain the same for all values of $t$, up to smooth isotopies inside the
reference fiber. For this reason, we do not need to make a
specific choice of $\lambda_i$. To fix ideas, let us say that $\lambda_0$
is close to $3$, $\lambda_1$ is close to $3j^2$, $\lambda_2$ is close to
$3j$, and $\lambda_i$ is close to infinity for $i\ge 3$; we again
choose the origin as base point, and note that the smooth elliptic curve
$W_k^{-1}(0)$ is diffeomorphic to $\ol{W_0}{}^{-1}(0)$, so we implicitly
identify them and again call our reference fiber $\Sigma_0$. We also
choose an ordered collection of arcs $\gamma_i$ joining the origin to
$\lambda_i$ which lie close to those considered in \S \ref{ss:critw0},
thus ensuring that the homology classes $[L_0],\dots,[L_{k+2}]\in
H_1(\Sigma_0,\Z)$ of the corresponding vanishing cycles remain those given
by Lemma \ref{l:vch1}.

Fixing a symplectic form $\omega_k$ on $M_k$ (compatible with $W_k$,
i.e.\ restricting positively to the fibers),
the vanishing cycles $L_0,\dots,L_{k+2}$
associated to the arcs $\gamma_0,\dots,\gamma_{k+2}$ naturally become
{\it Lagrangian submanifolds} of the reference fiber
$(\Sigma_0,\omega_{k|\Sigma_0})$ (cf.\ e.g.\ \cite{Ar,Se1,SeBook}).
Indeed, the symplectic form defines a natural horizontal distribution
outside of the critical points of $W_k$, given by the symplectic
orthogonal to the fiber.
Using this horizontal distribution, parallel transport induces
symplectomorphisms between the smooth fibers, and the vanishing cycle $L_i$
is by definition the set of points in the reference fiber $\Sigma_0$ for
which parallel transport along $\gamma_i$ converges to the critical point
in the fiber $W_k^{-1}(\lambda_i)$. It is also useful to consider the {\it
Lefschetz thimble} $D_i$, which is the set of points swept out by parallel
transport of $L_i$ above $\gamma_i$; by construction, $D_i$ is a
Lagrangian disk in $(M_k,\omega_k)$, fibered above the arc $\gamma_i$,
and $\partial D_i=L_i$.

We recall the following classical result (we provide a proof for
completeness):

\begin{lemma}\label{l:exactdeform}
A deformation of the system of arcs $\{\gamma_i\}$ by an isotopy in
$\mathrm{Diff}(\C,\mathrm{crit}(W_k))$ affects the vanishing cycles $L_i$
by Hamiltonian isotopies; moreover, the same property holds if the
symplectic fibration $(M_k,\omega_k,W_k)$ is deformed in a manner such that
the cohomology class $[\omega_k]$ remains constant and the critical points
of $W_k$ remain isolated and non-degenerate.
\end{lemma}

\proof 
We first consider a deformation of the system of arcs $\{\gamma_i\}$, based
at a regular value $\lambda_*\in\C\setminus\mathrm{crit}(W_k)$ (in our
case the origin), to an isotopic system of arcs $\{\gamma'_i\}$ based at a regular value
$\lambda'_*$. This means that we are given an arc $\delta:[0,1]\to
\C\setminus\mathrm{crit}(W_k)$ joining $\lambda_*$ to $\lambda'_*$,
and continuous families of arcs $\{\gamma_{i,t}\}$, $0\le t\le 1$, with
$\gamma_{i,0}=\gamma_i$ and $\gamma_{i,1}=\gamma'_i$, such that
$\gamma_{i,t}$ joins the regular value $\delta(t)$ to the critical value
$\lambda_i$, and
$\{\gamma_{i,t}\}_{0\le i\le k+2}$ is an ordered collection of arcs
for all $t\in [0,1]$.
The vanishing cycles $L'_i$ associated to the arcs $\gamma'_i$ live inside
$\Sigma'_*=W_k^{-1}(\lambda'_*)$, while the original vanishing cycles $L_i$
associated to $\gamma_i$ are submanifolds of $\Sigma_*=W_k^{-1}(\lambda_*)$.
However, we claim that the isotopy induces a symplectomorphism
$\phi:\Sigma_*\to \Sigma'_*$ with the property that $\phi(L_i)$ and $L'_i$
are mutually Hamiltonian isotopic for all $i$; this is the meaning of the
statement of the lemma.

Namely, parallel transport along the arc $\delta$ (using the
horizontal distribution described above) induces a symplectomorphism
$\phi$ from $\Sigma_*=W_k^{-1}(\lambda_*)$ to
$\Sigma'_*=W_k^{-1}(\lambda'_*)$. For all $t\in [0,1]$ we can consider the
vanishing cycle $L_{i,t}\subset W_k^{-1}(\delta(t))$ associated to the arc
$\gamma_{i,t}$, and its image $L'_{i,t}\subset \Sigma'_*$ under the
symplectomorphism induced by parallel transport along $\delta([t,1])$.
The family $L'_{i,t}$, $t\in [0,1]$ defines a smooth isotopy from
$L'_{i,0}=\phi(L_i)$ to $L'_{i,1}=L'_i$ through Lagrangian submanifolds
of $\Sigma'_*$. Moreover, each vanishing cycle $L_{i,t}\subset W_k^{-1}(
\delta(t))$ bounds a Lagrangian thimble $D_{i,t}$, and the cylinder
$C_{i,t}$ swept by $L_{i,t}$ under parallel transport along $\delta([t,1])$
is also Lagrangian. By continuity, the relative cycles $D_{i,t}\cup C_{i,t}$
(with boundary $L'_{i,t}$) all represent the same relative homotopy class
in $\pi_2(M_k,\Sigma'_*)$, and since $D_{i,t}$ and $C_{i,t}$ are Lagrangian
they all have zero symplectic area. This implies that the Lagrangian
submanifolds $L'_{i,t}$, $t\in [0,1]$ are all Hamiltonian isotopic inside
$\Sigma'_*$; in particular, $\phi(L_i)$ and $L'_i$ are Hamiltonian isotopic.

We now consider a symplectic fibration $W'_k:(M_k,\omega'_k)\to \C$ which
is isotopic to $W_k$ through an isotopy $W_{k,t}:(M_k,\omega_{k,t})\to\C$
that preserves the cohomology class of the symplectic form (i.e.,
$[\omega_{k,t}]=[\omega_k]$ for all $t\in [0,1]$). We assume that each
$W_{k,t}$ has isolated non-degenerate critical points. This allows us to 
deform the system of arcs $\{\gamma_i\}$ through a family
$\{\gamma_{i,t}\}$ with end points at the critical values of $W_{k,t}$;
for $t=1$ we obtain a system of arcs $\{\gamma'_i\}$ based at a regular
value $\lambda'_*$ of $W'_k$.
By Moser's theorem, there exists a continuous family of symplectomorphisms
$\phi_t$ from $(M_k,\omega_{k,t})$ to $(M_k,\omega'_k)$, or rather, since
these are non-compact manifolds, from open subsets of $(M_k,\omega_{k,t})$
to an open subset of $(M_k,\omega'_k)$; however,
after ``enlarging'' $(M_k,\omega'_k)$ by adding to $\omega'_k$ the pullback
of a suitable area form on $\C$, which affects neither the symplectic
structure on the fibers nor the parallel transport symplectomorphisms,
we can ensure that $\phi_t$ is defined over an arbitrarily large open subset
of $M_k$, which is good enough for our purposes. Moreover, by a relative version of
Moser's argument, we can also ensure that $\phi_t$ maps the reference fiber
of $W_{k,t}$ to the reference fiber of $W'_k$,
and in particular that $\phi=\phi_0$ maps
$\Sigma_*=W_k^{-1}(\lambda_*)$ to $\Sigma'_*={W'_k}^{-1}(\lambda'_*)$. 

We now claim that $\phi(L_i)\subset \Sigma'_*$ is Hamiltonian
isotopic to the vanishing cycle $L'_i$ of $W'_k$ associated to the arc
$\gamma'_i$. Indeed, by considering the images under $\phi_t$ of the
vanishing cycles $L_{i,t}$ associated to the arcs $\gamma_{i,t}$, we
obtain a smooth isotopy from $\phi(L_i)$ to $L'_i$ through
Lagrangian submanifolds of $\Sigma'_*$. Moreover, the thimbles $D'_i$ and
$\phi(D_i)$ represent the same relative homotopy class (as can be seen by
considering the images by $\phi_t$ of the thimbles $D_{i,t}$ associated to
$\gamma_{i,t}$), and both are
Lagrangian with respect to $\omega'_k$, which again implies that $\phi(L_i)$
and $L'_i$ are Hamiltonian isotopic.
\endproof

\subsection{A basis of $H_2(M_k)$}\label{ss:h2mk}


The manifold $M_k$ is simply connected, and its second Betti number is equal
to $k+2$. A $\mathbb{Q}$-basis of $H_2(M_k)$ is given by considering the homology
class of the fiber of $W_k$, $[\Sigma_0]$, and $k+1$ classes
$[\bar{C}],[\bar{C}_0],\dots,[\bar{C}_{k-1}]$ constructed from the vanishing
cycles $L_i$ and Lefschetz thimbles $D_i$ in the following manner.

By Lemma \ref{l:vch1} we have $[L_1]=[L_0]+[L_2]$ in $H_1(\Sigma_0,\Z)$,
so there exists
a 2-chain $C$ in $\Sigma_0$ such that $\partial C=-L_0+L_1-L_2$. Then
$$\bar{C}=C+D_0-D_1+D_2$$ is a $2$-cycle in $M_k$. Note that $[\bar{C}]$ is
in fact the image of the generator of $H_2((\C^*)^2,\Z)\simeq \Z$ under the inclusion map
(see the proof of Lemma 4.9 in \cite{AKO}). 

Similarly, for $0\le i<k$
we have $3\,[L_{3+i}]=[L_2]-[L_0]$ in $H_1(\Sigma_0,\Z)$, so there exists a 2-chain $C_i$ in
$\Sigma_0$ such that $\partial C_i=3\,L_{3+i}+L_0-L_2$, and we can consider
the $2$-cycle $$\bar{C}_i=C_i-3\,D_{3+i}-D_0+D_2$$ in $M_k$. We also introduce
2-chains $\Delta_{i,j}$ ($i,j\in\{0,\dots,k-1\}$) in $\Sigma_0$ such that
$\partial \Delta_{i,j}=L_{3+j}-L_{3+i}$, and the corresponding 2-cycles
$$\bar\Delta_{i,j}=\Delta_{i,j}+D_{3+i}-D_{3+j}.$$ We can choose $C_i$ and
$\Delta_{i,j}$ in such a way that $C_j-C_i=3\,\Delta_{i,j}$ (and hence
$[\bar{C}_j]-[\bar{C}_i]=3\,[\bar\Delta_{i,j}]$ in~$H_2(M_k)$).

To summarize the discussion, the vanishing cycles $L_i$ and the 2-chains
$C$, $C_i$, $\Delta_{i,j}$ are represented on Figures
\ref{fig:sigma0}--\ref{fig:deltaij}
(compare with Figure 2 in \cite{Se2} and with \cite{Ue}).

\begin{figure}[t]
\setlength{\unitlength}{0.9cm}
\begin{picture}(12,6.5)(0,-0.2)
\psset{unit=\unitlength,fillstyle=solid,linecolor=white}
\newrgbcolor{lt1}{0.9 0.9 0.9}
\newrgbcolor{lt2}{0.8 0.8 0.8}
\pspolygon[fillcolor=lt1](1.5,3)(3,2)(3,4)
\pspolygon[fillcolor=lt2](7.5,3)(6,2)(6,4)
\pspolygon[fillcolor=lt2](3,4)(3,6)(4.5,5)
\pspolygon[fillcolor=lt2](3,0)(3,2)(4.5,1)
\pspolygon[fillcolor=lt1](6,4)(6,6)(4.5,5)
\pspolygon[fillcolor=lt1](6,0)(6,2)(4.5,1)
\put(2.25,3){\small $-$}
\put(6.5,3){\small $+$}
\put(3.5,5){\small $+$}
\put(3.5,1){\small $+$}
\put(5.25,5){\small $-$}
\put(5.25,1){\small $-$}
\put(0,6){\line(1,-2){3}}
\put(3,0){\line(1,0){6}}
\put(0,6){\line(1,0){6}}
\put(6,6){\line(1,-2){3}}
\put(3.9,6){\vector(1,0){0}}
\put(4.1,6){\vector(1,0){0}}
\put(6.9,0){\vector(1,0){0}}
\put(7.1,0){\vector(1,0){0}}
\put(1.9,2.2){\vector(1,-2){0}}
\put(7.9,2.2){\vector(1,-2){0}}
\multiput(3,0)(0,0.2){30}{\line(0,1){0.1}}
\multiput(6,0)(0,0.2){30}{\line(0,1){0.1}}
\multiput(1.5,3)(0.66,0.44){7}{\line(3,2){0.45}}
\multiput(3,0)(0.66,0.44){7}{\line(3,2){0.45}}
\put(3,6){\line(3,-2){4.5}}
\put(1.5,3){\line(3,-2){4.5}}
\multiput(10.5,5)(0.8,0){2}{\line(1,0){0.54}}
\put(10.5,4){\line(1,0){1.35}}
\multiput(10.5,3)(0.2,0){7}{\line(1,0){0.1}}
\put(9.6,4.9){\large $L_0$}
\put(9.6,3.9){\large $L_1$}
\put(9.6,2.9){\large $L_2$}
\put(9.2,2){\large $L_{3+i}$}
\put(9.2,1.4){\large $L_{3+j}$}
\put(4.5,5){\circle*{0.1}}
\put(4.3,4.6){\small $z_0$}
\put(4.5,1){\circle*{0.1}}
\put(4.3,0.6){\small $y_0$}
\put(1.5,3){\circle*{0.1}}
\put(0.9,2.9){\small $x_0$}
\put(7.5,3){\circle*{0.1}}
\put(7.6,2.9){\small $x_0$}
\put(6,4){\circle*{0.1}}
\put(6.1,4.1){\small $y_1$}
\put(3,6){\circle*{0.1}}
\put(2.8,6.2){\small $x_1$}
\put(6,0){\circle*{0.1}}
\put(5.8,-0.4){\small $x_1$}
\put(3,2){\circle*{0.1}}
\put(3.1,2.1){\small $z_1$}
\put(6,6){\circle*{0.1}}
\put(6.1,6.1){\small $\bar{x}$}
\put(3,0){\circle*{0.1}}
\put(2.6,-0.3){\small $\bar{x}$}
\put(3,4){\circle*{0.1}}
\put(2.6,4.1){\small $\bar{y}$}
\put(6,2){\circle*{0.1}}
\put(6.1,1.7){\small $\bar{z}$}
\put(2.25,1.5){\line(1,2){2.25}}
\put(7.5,0){\line(1,2){0.75}}
\put(2.26,1.5){\line(1,2){2.25}}
\put(7.51,0){\line(1,2){0.75}}
\put(10.5,2){\line(1,0){1.35}}
\put(10.5,2.02){\line(1,0){1.35}}
\put(0.75,4.5){\line(1,2){0.75}}
\put(4.5,0){\line(1,2){2.25}}
\put(0.76,4.5){\line(1,2){0.75}}
\put(4.51,0){\line(1,2){2.25}}
\put(10.5,1.5){\line(1,0){1.35}}
\put(10.5,1.52){\line(1,0){1.35}}
\end{picture}
\caption{The vanishing cycles of $W_k$ and the chain $C$}\label{fig:sigma0}
\end{figure}
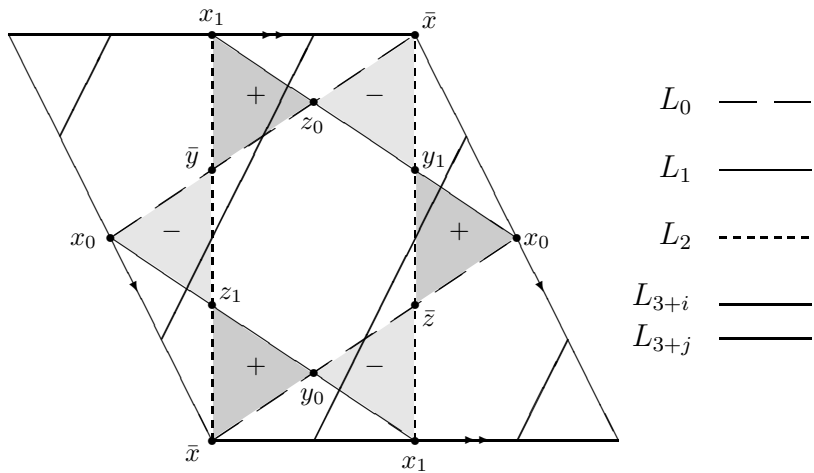

\begin{figure}[t]
\setlength{\unitlength}{0.9cm}
\begin{picture}(8.5,6.5)(0,0)
\psset{unit=\unitlength,fillstyle=solid,linecolor=white}
\newrgbcolor{lt1}{0.9 0.9 0.9}
\newrgbcolor{lt2}{0.8 0.8 0.8}
\pspolygon[fillcolor=lt1](3,0)(3,4)(6,6)(6,3)(5.25,1.5)
\pspolygon[fillcolor=lt2](6,3)(5.25,1.5)(6,2)
\pspolygon[fillcolor=lt1](6,3)(6.75,4.5)(7.5,3)(6,2)
\pspolygon[fillcolor=lt1](5.25,1.5)(6,2)(6,0)(4.5,0)
\pspolygon[fillcolor=lt2](3,0)(5.25,1.5)(4.5,0)
\pspolygon[fillcolor=lt2](6,6)(6.75,4.5)(6,3)
\pspolygon[fillcolor=lt2](0,6)(0.75,4.5)(1.5,6)
\pspolygon[fillcolor=lt1](1.5,6)(3,6)(3,4)(1.5,3)(0.75,4.5)
\put(4.4,3){\small $+$}
\put(5.65,2.05){\small -2}
\put(6.5,3){\small $-$}
\put(5.4,0.7){\small $-$}
\put(4,0.4){\small $+2$}
\put(6.1,4.4){\small $+2$}
\put(0.5,5.3){\small $+2$}
\put(1.8,4.5){\small $-$}
\put(0,6){\line(1,-2){3}}
\put(3,0){\line(1,0){6}}
\put(0,6){\line(1,0){6}}
\put(6,6){\line(1,-2){3}}
\put(3.9,6){\vector(1,0){0}}
\put(4.1,6){\vector(1,0){0}}
\put(6.9,0){\vector(1,0){0}}
\put(7.1,0){\vector(1,0){0}}
\put(1.9,2.2){\vector(1,-2){0}}
\put(7.9,2.2){\vector(1,-2){0}}
\multiput(3,0)(0,0.2){30}{\line(0,1){0.1}}
\multiput(6,0)(0,0.2){30}{\line(0,1){0.1}}
\multiput(1.5,3)(0.66,0.44){7}{\line(3,2){0.45}}
\multiput(3,0)(0.66,0.44){7}{\line(3,2){0.45}}
\put(3,6){\line(3,-2){4.5}}
\put(1.5,3){\line(3,-2){4.5}}
\put(6,6){\circle*{0.1}}
\put(6.1,6.1){\small $\bar{x}$}
\put(3,0){\circle*{0.1}}
\put(2.6,-0.3){\small $\bar{x}$}
\put(3,4){\circle*{0.1}}
\put(2.6,4.1){\small $\bar{y}$}
\put(6,2){\circle*{0.1}}
\put(6.1,1.7){\small $\bar{z}$}
%
\put(0.75,4.5){\line(1,2){0.75}}
\put(4.5,0){\line(1,2){2.25}}
\put(0.76,4.5){\line(1,2){0.75}}
\put(4.51,0){\line(1,2){2.25}}
\put(4.875,0.75){\circle*{0.1}}
\put(4.82,0.35){\small $b_i$}
\put(5.25,1.5){\circle*{0.1}}
\put(4.85,1.6){\small $a_i$}
\put(6,3){\circle*{0.1}}
\put(5.6,3.05){\small $c_i$}
\put(6.375,3.75){\circle*{0.1}}
\put(6.53,3.72){\small $b'_i$}
\end{picture}
\hspace{-1cm}
\begin{picture}(8.5,6.5)(0.5,0)
\psset{unit=\unitlength,fillstyle=solid,linecolor=white}
\newrgbcolor{lt1}{0.9 0.9 0.9}
\newrgbcolor{lt2}{0.8 0.8 0.8}
\pspolygon[fillcolor=lt1](0.75,4.5)(1.5,6)(4.5,6)(2.25,1.5)
\pspolygon[fillcolor=lt1](4.5,0)(6.75,4.5)(8.25,1.5)(7.5,0)
\put(0,6){\line(1,-2){3}}
\put(3,0){\line(1,0){6}}
\put(0,6){\line(1,0){6}}
\put(6,6){\line(1,-2){3}}
\put(3.9,6){\vector(1,0){0}}
\put(4.1,6){\vector(1,0){0}}
\put(6.9,0){\vector(1,0){0}}
\put(7.1,0){\vector(1,0){0}}
\put(1.9,2.2){\vector(1,-2){0}}
\put(7.9,2.2){\vector(1,-2){0}}
\multiput(3,0)(0,0.2){30}{\line(0,1){0.1}}
\multiput(6,0)(0,0.2){30}{\line(0,1){0.1}}
\multiput(1.5,3)(0.66,0.44){7}{\line(3,2){0.45}}
\multiput(3,0)(0.66,0.44){7}{\line(3,2){0.45}}
\put(3,6){\line(3,-2){4.5}}
\put(1.5,3){\line(3,-2){4.5}}
%
\put(2.25,1.5){\line(1,2){2.25}}
\put(7.5,0){\line(1,2){0.75}}
\put(2.26,1.5){\line(1,2){2.25}}
\put(7.51,0){\line(1,2){0.75}}
\put(0.75,4.5){\line(1,2){0.75}}
\put(4.5,0){\line(1,2){2.25}}
\put(0.76,4.5){\line(1,2){0.75}}
\put(4.51,0){\line(1,2){2.25}}
\end{picture}
\caption{The chains $C_i$ (left) and $\Delta_{i,j}$ (right)}\label{fig:deltaij}
\end{figure}
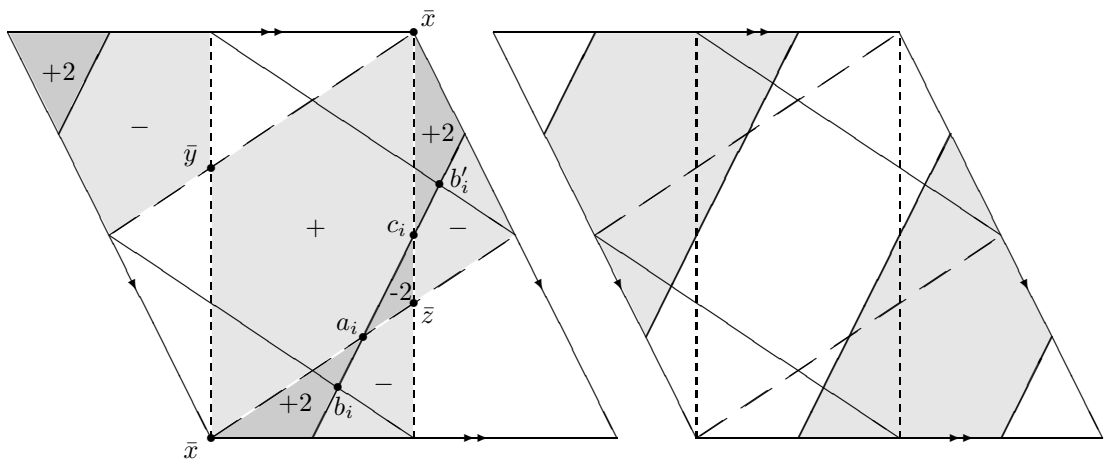


\section{Categories of vanishing cycles}\label{s:lagvc}

\subsection{Definition}\label{ss:deflagvc}

As proposed by Kontsevich \cite{KoEns} and Hori-Iqbal-Vafa \cite{HIV},
the category of A-branes associated to a Landau-Ginzburg model
$W:(M,\omega)\to\C$ is a Fukaya-type category which contains not only
compact Lagrangian submanifolds of $M$ but also certain non-compact
Lagrangians whose ends fiber in a specific way above half-lines in $\C$.
In the case where the critical points of $W$ are isolated and
non-degenerate, this category admits an exceptional collection whose
objects are Lagrangian thimbles associated to the critical points.
Following the formalism introduced by Seidel \cite{Se1,SeBook}, we view
it as the derived category of a finite directed
$A_\infty$-category $\FS(W,\{\gamma_i\})$ associated to an ordered
collection of arcs $\{\gamma_i\}$. We briefly recall the definition;
the reader is referred to \cite{Se1,SeBook} and to \S 3.1 of \cite{AKO}
for details.

Consider a symplectic fibration $W:(M,\omega)\to \C$ with isolated
non-degenerate critical points, and assume for simplicity that the
critical values $\lambda_0,\dots,\lambda_r$ of $W$ are distinct.
Pick a regular value $\lambda_*$ of $W$, and choose a collection of
arcs $\gamma_0,\dots,\gamma_r\subset\C$ joining $\lambda_*$ to the
various critical values of $W$, intersecting each other only at $\lambda_*$,
and ordered in the clockwise direction around $\lambda_*$. Consider the
horizontal distribution defined by the symplectic form: by parallel
transport along the arc $\gamma_i$, we obtain a Lagrangian thimble $D_i$
and a vanishing cycle $L_i=\partial D_i\subset \Sigma_*$ (where
$\Sigma_*=W^{-1}(\lambda_*)$). After a small perturbation we can
always assume that the vanishing cycles $L_i$ intersect each other
transversely inside $\Sigma_*$.
The following definition is motivated by the observation that
the intersection theory of the Lagrangian thimbles
$D_i\subset M$ is closely related to that of the vanishing cycles $L_i$
inside $\Sigma_*$ \cite{Se1}:

\begin{defi}[Seidel]\label{def:fs}
The directed category of vanishing cycles $\FS(W,\{\gamma_i\})$ is an
$A_\infty$-category (over a coefficient ring $R$) with objects
$L_0,\dots,L_r$ corresponding to the vanishing cycles (or more accurately
to the thimbles); the morphisms between
the objects are given by
$$\Hom(L_i,L_j)=\begin{cases}
CF^*(L_i,L_j;R)=R^{|L_i\cap L_j|} & \mathrm{if}\ i<j\\
R\cdot id & \mathrm{if}\ i=j\\
0 & \mathrm{if}\ i>j;
\end{cases}$$
and the differential $m_1$, composition $m_2$ and higher order
products $m_k$ are defined in terms of Lagrangian Floer homology inside
$\Sigma_*$. More precisely,
$$m_k:\Hom(L_{i_0},L_{i_1})\otimes \dots\otimes
\Hom(L_{i_{k-1}},L_{i_k}) \to \Hom(L_{i_0},L_{i_k})[2-k]$$
is trivial when the inequality $i_0<i_1<\dots<i_k$ fails to hold (i.e.\ it
is always zero in this case, except for $m_2$ where composition with an
identity morphism is given by the obvious formula).
When $i_0<\dots<i_k$, $m_k$ is defined by fixing a generic
$\omega$-compatible almost-complex structure on $\Sigma_*$ and counting
pseudo-holomorphic maps from a disk with $k+1$ cyclically ordered
marked points on its boundary to $\Sigma_*$, mapping the marked points
to the given intersection points between vanishing cycles, and the portions
of boundary between them to $L_{i_0},\dots,L_{i_k}$ respectively.
\end{defi}

This definition calls for several clarifications. First of all, in
our case $\Sigma_*$ is a smooth elliptic curve and the vanishing cycles
are homotopically non-trivial closed loops, we have $\pi_2(\Sigma_*)=0$ and
$\pi_2(\Sigma_*,L_i)=0$; hence, we need not be concerned by bubbling issues
in the definition of the Floer differential and products. In fact, the
pseudo-holomorphic disks in $\Sigma_*$ that we have to consider are nothing
but immersed polygonal regions bounded by the vanishing cycles, satisfying
a local convexity condition at each corner point.

Also, the Maslov class vanishes identically, so we have a well-defined
$\Z$-grading by Maslov index on the Floer complexes $CF^*(L_i,L_j;R)$
once we choose graded Lagrangian lifts of the vanishing cycles.
Since in our case $c_1(\Sigma_*)=0$, we can do this by considering a nowhere
vanishing  1-form $\Omega\in\Omega^1(\Sigma_*,\C)$ and choosing a real
lift of the phase function $\phi_i=\mathrm{arg}(\Omega_{|L_i}):L_i\to S^1$
for each vanishing cycle. The degree of
a given intersection point $p\in L_i\cap L_j$ is then determined by
the difference between the phases of $L_i$ and $L_j$ at $p.$

Our next remark is that the pseudo-holomorphic disks appearing in
Definition \ref{def:fs} are counted with appropriate weights, and with
signs determined by choices of orientations of the relevant moduli spaces.
The orientation is determined by the choice of a spin structure for each
vanishing cycle $L_i$; in our case this spin structure must extend to
the thimble, so it is necessarily the non-trivial one. In the
one-dimensional case there is a convenient recipe for determining the
correct sign factors, due to Seidel \cite{SeBook}. As will be clear from
the discussion in \S \ref{ss:objmor} below, we will only be interested
in the specific case where all morphisms have even degree and all spin
structures are non-trivial. The sign rule can then be summarized as
follows: pick a marked point on each $L_i$, distinct from the intersections
with the other vanishing cycles; then the sign associated to a
pseudo-holomorphic map $u:(D^2,\partial D^2)\to (\Sigma_*,\cup L_i)$ is
$(-1)^{\nu(u)}$, where $\nu(u)$ is the number of marked points that the
boundary of $u$ passes through (\cite{SeBook}, see also \S 4.6 of \cite{AKO}).

Finally, the weight attributed to each pseudo-holomorphic map $u$ keeps track
of its relative homology class, which makes it possible to avoid convergence
problems. The usual approach favored by mathematicians is to work over a
Novikov ring, which keeps track of the relative homology class by introducing
suitable formal variables. To remain closer to the physics, we use $\C$ as
our coefficient ring, and assign weights according to the symplectic areas;
this is in fact equivalent to working over the Novikov ring and specializing
at the cohomology class of the symplectic form.

The weight formula is simplest when there is no B-field; in that case, we
consider untwisted Floer theory, since any flat unitary bundle over the
thimble $D_i$ is trivial and hence restricts to $L_i$ as the trivial bundle.
We then count each map $u:(D^2,\partial D^2)\to (\Sigma_*,\cup L_i)$ with a
coefficient $(-1)^{\nu(u)} \exp(-2\pi\int_{D^2} u^*\omega)$. (The
normalization factor $2\pi$ is purely a matter of conventions, and is
sometimes omitted in the literature; here we include it for convenience).
Hence, given two intersection points $p\in L_i\cap L_j$, $q\in L_j\cap
L_k$ ($i<j<k$), we have by definition
$$m_2(p,q)=\sum_{\substack{r\in L_i\cap L_k\\\deg r=\deg p+\deg q}}
\Biggl(\sum_{[u]\in \mathcal{M}(p,q,r)}\!\!(-1)^{\nu(u)}
\exp(-2\pi\int_{D^2} u^*\omega) \Biggr)\,r$$
where $\mathcal{M}(p,q,r)$ is the moduli space of pseudo-holomorphic maps
$u$ from the unit disk to $\Sigma_*$ (equipped with a generic $\omega$-compatible
almost-complex structure) such that $u(1)=p$, $u(\mathrm{j})=q$,
$u(\mathrm{j}^2)=r$ (where $\mathrm{j}=\exp(\frac{2i\pi}{3})$), and mapping
the portions of unit circle $[1,\mathrm{j}]$, $[\mathrm{j},\mathrm{j}^2]$,
$[\mathrm{j}^2,1]$ to $L_i$, $L_j$ and $L_k$ respectively.
The other products are defined similarly. (Observe that Seidel's definition
\cite{Se1} does not involve any weights; this is because he only considers
exact Lagrangian submanifolds in exact symplectic manifolds, in
which case the symplectic areas are entirely determined by the primitives
of the Liouville form and can be eliminated by considering suitably rescaled
bases of the Floer complexes.)

In presence of a B-field, the weights are modified by the fact that we now
consider {\it twisted} Floer homology. Indeed, each thimble $D_i$ now comes
equipped with a trivial complex line bundle $E_i=\underline{\C}$
and a connection $\nabla_i$ with curvature $-2\pi iB$, so its boundary $L_i$ is
equipped with the restricted bundle and the restricted connection,
whose holonomy is $\mathrm{hol}_{\nabla_i}(L_i)=\exp(-2\pi i\int_{D_i} B)$
by Stokes' theorem. Since this property characterizes the connection
$\nabla_i$ uniquely up to gauge, we can drop the line
bundle and the connection from the notation when considering the objects
$(L_i,E_i,\nabla_i)$ of $\FS(W,\{\gamma_i\})$. Nonetheless, the holonomy
of $\nabla_i$ modifies the weights attributed to the pseudo-holomorphic
disks in the definition of the twisted Floer differentials and compositions.
Namely, the weight attributed to a given
pseudo-holomorphic map $u:(D^2,\partial D^2)\to (\Sigma_*,\cup L_i)$
is modified by a factor corresponding to the holonomy along its boundary,
and becomes
$$(-1)^{\nu(u)}\,\mathrm{hol}(u(\partial D^2))\,\exp(2\pi i\int_{D^2}
u^*(B+i\omega)).$$
More precisely, we fix trivializations of the line bundles $E_i$, so that
for each intersection point $p\in L_i\cap L_j$ we have a preferred
isomorphism between the fibers $(E_i)_{|p}$ and $(E_j)_{|p}$;
then it becomes possible to define the holonomy along the closed loop
$u(\partial D^2)$ using the parallel transport induced by $\nabla_i$ from
one ``corner'' of $u$ to the next one, and the chosen isomorphism at each
corner.

Although the category $\FS(W,\{\gamma_i\})$ depends on the chosen ordered
collection of arcs $\{\gamma_i\}$, Seidel has obtained the following result
\cite{Se1} (for the exact case, but the proof extends to our situation):

\begin{theo}[Seidel]
If the ordered collection $\{\gamma_i\}$ is
replaced by another one $\{\gamma'_i\}$, then the categories
$\FS(W,\{\gamma_i\})$ and $\FS(W,\{\gamma'_i\})$ differ by a sequence of
mutations.
\end{theo}

Hence, the category naturally associated to the
fibration $W$ is not the finite $A_\infty$-category defined above,
but rather a (bounded) {\it derived} category, obtained from
$\FS(W,\{\gamma_i\})$
by considering twisted complexes of formal direct sums of Lagrangian
vanishing cycles, and adding idempotent splittings and formal inverses
of quasi-isomorphisms (see \cite{KoEns} and \S 5 of \cite{Se1}). 
If two categories differ by mutations, then their derived categories are
equivalent; hence the derived category $\db{\FS(W)}$ depends only on the
symplectic topology of $W$ and not on the choice of an ordered system
of arcs \cite{Se1}. 

For the examples we consider, the $A_\infty$-category $\FS(W,\{\gamma_i\})$
will in fact be an honest category (see below); the bounded derived category
$\db{\FS(W)}$ is then by definition the bounded derived category
of finite rank modules over the algebra associated to this category.

\subsection{Objects and morphisms}\label{ss:objmor}

We now determine the categories $\FS(W_k,\{\gamma_i\})$ associated to the
Landau-Ginzburg models $(M_k,W_k)$ mirror to Del Pezzo surfaces and the
systems of arcs $\{\gamma_i\}$ introduced in \S \ref{ss:critwk}. We start
with the objects and morphisms.

Recall that $W_k$ has $k+3$ isolated critical points, giving rise to $k+3$
vanishing cycles $L_0,\dots,L_{k+2}$ in the reference fiber $\Sigma_0\simeq
W_k^{-1}(0)$. The homology classes of these vanishing cycles have been
determined in \S \ref{s:mirrors} and are given by Lemma \ref{l:vch1}; these
determine the vanishing cycles up to Lagrangian isotopy.

The derived category of vanishing cycles is not affected if we modify some
of the vanishing cycles by Hamiltonian isotopies (more precisely, a
Hamiltonian isotopy induces a chain map on the Floer complexes, which
yields a quasi-isomorphism between the finite
$A_\infty$-categories of vanishing cycles). Hence, equipping the elliptic
curve $\Sigma_0$ with a compatible flat metric, we can identify $\Sigma_0$
with the quotient of $\C$ by a lattice, and represent the vanishing cycles
$L_i$ by closed geodesics parallel to those represented in Figure
\ref{fig:sigma0}.

Assume that the cohomology class of the symplectic form $\omega_k$ on
$M_k$ is generic (or more precisely, with the notations of \S \ref{ss:h2mk},
that $[\omega_k]\cdot[\bar\Delta_{i,j}]$ is never an integer multiple of
$[\omega_k]\cdot[\Sigma_0]$). Then the geodesics $L_i$ are all distinct, and their
intersections are as pictured in Figure \ref{fig:sigma0}, so we have:

\begin{lemma} The geometric intersection numbers between the vanishing
cycles are:

$\bullet$ $|L_0\cap L_1|=|L_0\cap L_2|=|L_1\cap L_2|=3$;

$\bullet$ for $0\le i<k$, $|L_0\cap
L_{3+i}|=|L_2\cap L_{3+i}|=1$ and $|L_1\cap L_{3+i}|=2$;

$\bullet$ for $0\le i<j<k$, $|L_{3+i}\cap L_{3+j}|=0$ as soon as
$[\omega_k]\cdot[\Sigma_0] \not| \ \,[\omega_k]\cdot[\bar\Delta_{i,j}] $.
\end{lemma}

In the rest of this section, unless otherwise specified we always assume
that the vanishing cycles $L_i$ are represented by distinct closed geodesics.

As in \cite{AKO}, we denote by $x_0,y_0,z_0$ (resp.\ $x_1,y_1,z_1$ and
$\bar{x},\bar{y}, \bar{z}$) the generators of $\Hom(L_0,L_1)$ (resp.\
$\Hom(L_1,L_2)$ and $\Hom(L_0,L_2)$) corresponding to the intersection
points represented in Figure \ref{fig:sigma0}. Moreover, we denote by
$a_i$ (resp.\ $b_i,b'_i$ and $c_i$) the generators of $\Hom(L_0,L_{3+i})$
(resp.\ $\Hom(L_1,L_{3+i})$ and $\Hom(L_2,L_{3+i})$) corresponding to the
intersection points between these vanishing cycles (see Figure
\ref{fig:deltaij}).

\begin{lemma}
For suitable choices of graded lifts of the vanishing cycles, all the
morphisms in $\FS(W_k,\{\gamma_i\})$ have degree $0$.
\end{lemma}

\proof
Equip $\Sigma_0$ with a compatible flat metric and with a constant
holomorphic 1-form $\Omega$. Taking geodesic representatives of the
vanishing cycles, the phase functions $\phi_i=\arg(\Omega_{|L_i}):L_i\to
\R/2\pi\Z$ are constant, and we can normalize $\Omega$ so that it takes real
negative values on the oriented tangent space to $L_0$, i.e.\ $\phi_0=\pi$.
Then it is possible to choose real lifts $\tilde\phi_i\in\R$ of the phases
in such a way that $\pi=\tilde\phi_0>\tilde\phi_1>\tilde\phi_2>\tilde\phi_3
=\dots=\tilde\phi_{k+2}>0$ (see Figure \ref{fig:sigma0} and recall the
orientations chosen in Lemma \ref{l:vch1}). In the 1-dimensional case, the
relationship between Maslov index and phase is very simple: given a
transverse intersection point $p$ between two graded Lagrangians
$L,L'\subset \Sigma_0$, the Maslov index of $p\in CF^*(L,L')$ is equal to
the smallest integer greater than
$\frac{1}{\pi}(\tilde\phi_{L'}(p)-\tilde\phi_L(p))$. Since we only consider
the Floer complexes $CF^*(L_i,L_j)$ for $i<j$, which implies that
 $\tilde\phi_j-\tilde\phi_i\in (-\pi,0)$ at every intersection point,
for these choices of graded Lagrangian lifts of the vanishing cycles
all morphisms in $\FS(W_k,\{\gamma_i\})$ have degree 0.
\endproof

Since each product $m_j$ shifts degree by $2-j$, it follows immediately that
the $A_\infty$-category $\FS(W_k,\{\gamma_i\})$ is actually an honest
category:

\begin{cor}
$m_j=0$ for all $j\neq 2$.
\end{cor}

Hence, the final step of the argument is a careful study of the various immersed
triangular regions bounded by the vanishing cycles in $\Sigma_0$ and their
contributions to $m_2$.

\subsection{Compositions}\label{ss:m2}
As before we assume that the Lagrangian vanishing cycles are realized by
distinct closed geodesics in the flat torus $\Sigma_0$, and we determine
the contributions to $m_2$ of the various immersed triangular regions in
$(\Sigma_0,\cup L_i)$. We use the notations introduced in \S \ref{ss:objmor}
for the intersection points, and those introduced in \S \ref{ss:h2mk} for
various 2-chains in $\Sigma_0$ and the corresponding 2-cycles in $M_k$.
We also introduce the following notations:

\begin{defi}\label{def:zetas} Let 
$q_C=\exp(2\pi i\,[B+i\omega]\cdot[\bar{C}])$ and
$q_F=\exp(2\pi i\,[B+i\omega]\cdot[\Sigma_0])$, and define
$$\zeta_+=\sum_{n\in\Z}(-1)^n\,q_C^{n\vphantom{/}}\,q_F^{n(3n+1)/2},
\ \ \zeta_-=\sum_{n\in\Z}(-1)^n\,q_C^{n\vphantom{/}}\,q_F^{n(3n-1)/2},
\ \ \zeta_0=\sum_{n\in\Z}(-1)^n\,q_C^{n\vphantom{/}}\,q_F^{3n(n-1)/2}.$$
\end{defi}

\noindent
Since $\omega$ is a symplectic form on $\Sigma_0$, we have
$|q_F|=\exp(-2\pi\,[\omega]\cdot[\Sigma_0])<1$, which ensures the convergence
of the series $\zeta_+$, $\zeta_-$ and $\zeta_0$. 

\begin{prop}\label{prop:m2p2}
There exist constants $\alpha_{xy},\alpha_{yx},\alpha_{yz},\alpha_{zy},
\alpha_{zx},\alpha_{xz}\in\C$ such that
\begin{eqnarray*}
m_2(x_0,y_1)=\alpha_{xy}\bar{z},&\quad& m_2(y_0,x_1)=\alpha_{yx}\bar{z},\\
m_2(y_0,z_1)=\alpha_{yz}\bar{x},&\quad& m_2(z_0,y_1)=\alpha_{zy}\bar{x},\\
m_2(z_0,x_1)=\alpha_{zx}\bar{y},& & m_2(x_0,z_1)=\alpha_{xz}\bar{y},
\end{eqnarray*}
and these constants satisfy the relation
\begin{equation}\label{eq:alpha46}\displaystyle
\frac{\alpha_{xy}\alpha_{yz}\alpha_{zx}}{\alpha_{yx}\alpha_{zy}
\alpha_{xz}}=-q_C\left(\frac{\sum_{n\in\Z}\, (-1)^n \,q_C^{n\vphantom{/}}\,
q_F^{n(3n+1)/2}
} {\sum_{n\in\Z}\, (-1)^n\, q_C^{n\vphantom{/}}\, q_F^{n(3n-1)/2}
}\right)^{\!3}=-q_C\left(\frac{\zeta_+}{\zeta_-}\right)^{\!3}.
\end{equation}
\end{prop}

\begin{remark}
The quantity appearing in the right-hand side of (\ref{eq:alpha46})
can be understood in terms of certain theta functions; see \S \ref{ss:theta}
for details.
\end{remark}

Before giving the proof, we make an observation which will be useful
throughout this section. The geodesics $L_i$ are not necessarily those
pictured in Figure \ref{fig:sigma0}, but they are parallel to them. So
we can deform (in a non-Hamiltonian manner) the configuration of vanishing
cycles to that of Figure \ref{fig:sigma0}, and all intersection points
and relative 2-cycles in $(\Sigma_0,\cup L_i)$ can be
followed through the deformation. Hence, immersed triangular regions
in $(\Sigma_0,\cup L_i)$ are in one to one correspondence with those in
the configuration of Figure \ref{fig:sigma0} (but of course the deformation
does not preserve areas). Moreover, we can choose the deformation in such
a way that the relative 2-cycles $C$ and $C_i$ in $\Sigma_0$ deform to those
represented on Figures \ref{fig:sigma0}--\ref{fig:deltaij} (rather than to
2-cycles which differ by a multiple of the fundamental cycle of $\Sigma_0$).

\proof[Proof of Proposition \ref{prop:m2p2}]
The composition $m_2(x_0,y_1)$ is the sum of an infinite series of
contributions, corresponding to all immersed triangular regions in
$\Sigma_0$ with corners at the intersection points $x_0$, $y_1$, and one
of the points in $L_0\cap L_2$. By the above remark we can enumerate these
regions by looking at Figure \ref{fig:sigma0}. Considering the side which lies
on $L_1$, it is then easy to see that for every homotopy class of arc joining
$x_0$ to $y_1$ inside $L_1$ there is a unique such immersed triangular
region, and the third vertex is always $\bar{z}$.

These various regions can be labelled by integers $n\in\Z$ in such a way
that, denoting by $T_{xy,n}$ the corresponding
2-chains in $\Sigma_0$, we have $\partial T_{xy,n}-\partial
T_{xy,n'}=(n-n')(-L_0+L_1-L_2)$ for all $n,n'\in\Z$.
We can choose the integer labels in such a way that, after deforming to
the configuration in Figure \ref{fig:sigma0}, $T_{xy,0}$ becomes
the smallest triangle with vertices $x_0,y_1,\bar{z}$. (So, in Figure
\ref{fig:sigma0}, $T_{xy,-1}$ is the immersed region bounded by the
portions of $L_0\cup L_1\cup L_2$ which do not belong to $\partial T_{xy,0}$;
and all the other $T_{xy,n}$ have edges which wrap more than once around
the vanishing cycles).

By comparing $\partial T_{xy,n}$ and $\partial T_{xy,0}$, it is clear that
the 2-chain represented by $T_{xy,n}$ can be expressed in the form
$T_{xy,n}=T_{xy,0}+nC+\phi(n)\Sigma_0$ for some $\phi(n)\in\Z$. Moreover, by
looking at Figure \ref{fig:sigma0} one easily checks that
$\phi(n)=\frac{1}{2}n(3n+1)$. (So e.g.\ $T_{xy,-1}=T_{xy,0}-C+\Sigma_0$, and
$T_{xy,1}=T_{xy,0}+C+2\Sigma_0$). Let $\psi_{xy}\in\C$ be the coefficient of
the contribution of $T_{xy,0}$ to $m_2(x_0,y_1)$. Then, by comparing the
symplectic areas and boundary holonomies for $T_{xy,n}$ and $T_{xy,0}$,
one easily checks that the contribution of $T_{xy,n}$ is equal to
$$(-1)^n\exp\biggl(2\pi i\,[B+i\omega]\cdot
\Bigl(n[\bar{C}]+\frac{n(3n+1)}{2}[\Sigma_0]\Bigr)\biggr)\,\psi_{xy}=
(-1)^n\, q_C^{n\vphantom{/}}\,q_F^{n(3n+1)/2}\,\psi_{xy}.$$
In this expression the sign factor $(-1)^n$ is due to the non-triviality
of the spin structures (observe that $\partial C=-L_0+L_1-L_2$ passes once
through each of the three marked points on $L_0,L_1,L_2$); the total
holonomy of the flat connections $\nabla_i$ along $\partial
T_{xy,n}-\partial T_{xy,0}=n\,\partial C$ is
$\exp(2\pi i\,n\int_{D_0-D_1+D_2}B)$ by Stokes' theorem; and the integral of
$B+i\omega$ over $T_{xy,n}$ differs from that over $T_{xy,0}$ by
the amount $n\int_C (B+i\omega)+\frac{1}{2}n(3n+1)\,[B+i\omega]\cdot[\Sigma_0]$.

Summing over $n\in\Z$, and using the notation introduced in Definition
\ref{def:zetas}, we obtain
$$\alpha_{xy}=\zeta_+\,\psi_{xy}.$$
The calculations of $m_2(y_0,z_1)$ and $m_2(z_0,x_1)$ are exactly
identical, and lead to similar expressions. Namely, denote by
$\psi_{yz}$ (resp.\ $\psi_{zx}$) the contribution of the triangular
region $T_{yz,0}$ (resp.\ $T_{zx,0}$) which, after deforming to the
configuration in Figure \ref{fig:sigma0}, corresponds to the smallest
triangle with vertices $y_0,\,z_1,\,\bar{x}$ (resp.\ $z_0,\,x_1,\,\bar{y}$).
Then one easily checks by the same argument as above that
$\alpha_{yz}=\zeta_+\,\psi_{yz}$ and $\alpha_{zx}=\zeta_+\,\psi_{zx}$.

Next we consider the composition $m_2(y_0,x_1)$, which is again the sum of
an infinite series of contributions from triangular regions $T_{yx,n}$,
$n\in\Z$, which all have vertices $y_0,\,x_1,\,\bar{z}$. We can choose
the labels in such a way that, after deforming to the configuration in
Figure \ref{fig:sigma0}, $T_{yx,0}$ becomes the smallest such triangle,
and $T_{yx,n}=T_{yx,0}+nC+\frac{1}{2}n(3n-1)\Sigma_0$. Denoting by
$\psi_{yx}$ the coefficient associated to $T_{yx,0}$, it is easy to check
by the same argument as above that the contribution of $T_{yx,n}$ is equal to
$(-1)^n\,q_C^{n\vphantom{/}}\,q_F^{n(3n-1)/2}\,\psi_{yx}$, so that
$$\alpha_{yx}=\zeta_-\psi_{yx}.$$
Similarly, with the obvious notations we have
$\alpha_{zy}=\zeta_-\psi_{zy}$ and $\alpha_{xz}=\zeta_-\psi_{xz}$.
Finally, observe that
$$\frac{\psi_{xy}\psi_{yz}\psi_{zx}}{\psi_{yx}\psi_{zy}\psi_{xz}}
= -q_C.$$
Indeed, $T_{xy,0}+T_{yz,0}+T_{zx,0}-T_{yx,0}-T_{zy,0}-
T_{xz,0}=C$ (cf.\ Figure \ref{fig:sigma0}). Therefore, comparing the
weights associated to these various triangles, the weighting by area gives
a factor of $\exp(2\pi i\int_C B+i\omega)$, while the holonomy along
the boundary $\partial C=-L_0+L_1-L_2$ is equal to $\exp(2\pi i\int_{D_0-D_1+D_2}
B)$, and finally the minus sign is due to the orientation conventions, since
$\partial C$ passes once through each of the three marked points on the
vanishing cycles. Hence
$$\frac{\alpha_{xy}\alpha_{yz}\alpha_{zx}}{\alpha_{yx}\alpha_{zy}
\alpha_{xz}}=\frac{\psi_{xy}\psi_{yz}\psi_{zx}\,\zeta_+^3}{\psi_{yx}
\psi_{zy}\psi_{xz}\,\zeta_-^3}=-q_C\biggl(\frac{\zeta_+}{\zeta_-}\biggr)^3.$$
\endproof

\begin{remark}
If $[\omega+iB]\cdot[\bar{C}]=0$, then $q_C=1$ and the ratio
between $\alpha_{xy}\alpha_{yz}\alpha_{zx}$ and $\alpha_{yx}\alpha_{zy}
\alpha_{xz}$ becomes equal to $-1$ irrespective of the value of $q_F$;
this corresponds to a classical (commutative) Del Pezzo surface.

Moreover, in the limit where
$[\omega]\cdot[\Sigma_0]\to\infty$, we have $q_F=0$ and the ratio becomes
$-q_C$, which corresponds to the toric case studied in \cite{AKO}.

\end{remark}

\begin{prop}\label{prop:m2p2f}
There exist constants $\alpha_{xx},\alpha_{yy},\alpha_{zz}\in\C$ such that
$$m_2(x_0,x_1)=\alpha_{xx}\bar{x},\quad
m_2(y_0,y_1)=\alpha_{yy}\bar{y},\quad
m_2(z_0,z_1)=\alpha_{zz}\bar{z},$$
and these constants satisfy the relation
$$\displaystyle
\frac{\alpha_{xx}\alpha_{yy}\alpha_{zz}}{\alpha_{yx}\alpha_{zy}
\alpha_{xz}}=-\frac{q_F}{q_C}
\left(\frac{\sum_{n\in\Z}\, (-1)^n \,q_C^{n\vphantom{/}}\,
q_F^{3n(n-1)/2}
} {\sum_{n\in\Z}\, (-1)^n\, q_C^{n\vphantom{/}}\, q_F^{n(3n-1)/2}
}\right)^{\!3}=
-\frac{q_F}{q_C} \left(\frac{\zeta_0}{\zeta_-}\right)^{\!3}.$$
\end{prop}

\proof
The argument is similar to the proof of Proposition \ref{prop:m2p2}.
The immersed triangular regions which contribute to $m_2(x_0,x_1)$ all
have vertices $\bar{x}$ as their third vertex, and can be indexed by integers $n\in\Z$
in a manner such that $\partial T_{xx,n}-\partial T_{xx,n'}=(n-n')\,\partial
C$ for all $n,n'\in\Z$. We can choose the integer labels in such a way that,
after deforming to the standard configuration,
$T_{xx,0}$ and $T_{xx,1}=T_{xx,0}+C$ are the two embedded triangles with
vertices $x_0,\,x_1,\,\bar{x}$ visible on Figure \ref{fig:sigma0}.
It is then easy to check that
$T_{xx,n}=T_{xy,0}+nC+\frac{3}{2}n(n-1)\Sigma_0$.
Hence, denoting by $\psi_{xx}$ the coefficient associated to $T_{xx,0}$,
we have
$$\alpha_{xx}=\zeta_0\psi_{xx},$$
by the same argument as in previous calculations. Similarly, with the
obvious notations, we have $\alpha_{yy}=\zeta_0\psi_{yy}$ and
$\alpha_{zz}=\zeta_0\psi_{zz}$.
Moreover, $T_{xx,0}+T_{yy,0}+T_{zz,0}-T_{yx,0}-T_{zy,0}-T_{xz,0}=\Sigma_0-C$,
which implies (by the same argument as above) that
$$\frac{\psi_{xx}\psi_{yy}\psi_{zz}}{\psi_{yx}\psi_{zy}\psi_{xz}}
= -\frac{q_F}{q_C}.$$ Therefore
$$\frac{\alpha_{xx}\alpha_{yy}\alpha_{zz}}{\alpha_{yx}\alpha_{zy}
\alpha_{xz}}=\frac{\psi_{xx}\psi_{yy}\psi_{zz}\,\zeta_0^3}{\psi_{yx}
\psi_{zy}\psi_{xz}\,\zeta_-^3}=-\frac{q_F}{q_C}\biggl(\frac{\zeta_0}
{\zeta_-}\biggr)^3.$$
\endproof

When $q_F=0$ (in particular in the toric case) we have
$\alpha_{xx}\alpha_{yy}\alpha_{zz}=0$, as in \cite{AKO}. The same
conclusion also holds when $q_C=1$ (the commutative case). In fact,
when $q_C=1$ each of the constants $\alpha_{xx},\,\alpha_{yy},\,\alpha_{zz}$
is zero, since in that case we have $\zeta_0=0$ (because the terms
corresponding to $n$ and $1-n$ in the series defining $\zeta_0$
exactly cancel each other).

\begin{defi}\label{def:zetais} 
Let $q_i=\exp(2\pi i\,[B+i\omega]\cdot[\bar{C}_i])$, and define
$$\zeta_{i,+}=\sum_{n\in\Z}(-1)^n\,q_i^{n\vphantom{/}}\,q_F^{n(3n+1)/2},
\ \ \zeta_{i,-}=\sum_{n\in\Z}(-1)^n\,q_i^{n\vphantom{/}}\,q_F^{n(3n-1)/2},
\ \ \zeta_{i,0}=\sum_{n\in\Z}(-1)^n\,q_i^{n\vphantom{/}}\,q_F^{3n(n-1)/2}.$$
\end{defi}

\begin{prop}\label{prop:m2beta}
There exist constants $\beta_{\bar{x},i},\beta_{\bar{y},i},\beta_{\bar{z},i}
\in\C$ such that
$$m_2(\bar{x},c_i)=\beta_{\bar{x},i} a_i,\quad
m_2(\bar{y},c_i)=\beta_{\bar{y},i} a_i,\quad
m_2(\bar{z},c_i)=\beta_{\bar{z},i} a_i,$$
and these constants satisfy the relations
$$\frac{\beta_{\bar{z},i}^2\,\alpha_{xy}\alpha_{zz}}
{\beta_{\bar{x},i}\beta_{\bar{y},i}\,\alpha_{zy}\alpha_{xz}}
=\Bigl(\frac{\zeta_{i,-}}{\zeta_-}\Bigr)^2\frac{\zeta_+\,\zeta_0}
{\zeta_{i,+}\,\zeta_{i,0}},
$$
$$\frac{\beta_{\bar{x},i}^2\,\alpha_{yz}\alpha_{xx}}
{\beta_{\bar{y},i}\beta_{\bar{z},i}\,\alpha_{xz}\alpha_{yx}}
=-q_i\,\Bigl(\frac{\zeta_{i,+}}{\zeta_-}\Bigr)^2\frac{\zeta_+\,\zeta_0}
{\zeta_{i,0}\,\zeta_{i,-}},
\ \ \mathit{and}\ \
\frac{\beta_{\bar{y},i}^2\,\alpha_{zx}\alpha_{yy}}
{\beta_{\bar{z},i}\beta_{\bar{x},i}\,\alpha_{yx}\alpha_{zy}}
=-\frac{q_F}{q_i}\,\Bigl(\frac{\zeta_{i,0}}{\zeta_-}\Bigr)^2
\frac{\zeta_+\,\zeta_0}{\zeta_{i,-}\,\zeta_{i,+}},
$$
where 
$\zeta_+,\,\zeta_-,\,\zeta_0,\,\zeta_{i,+},\,\zeta_{i,-},\,\zeta_{i,0}$,
$q_i$ and $q_F$ are as in Definitions \ref{def:zetas} and \ref{def:zetais}.
\end{prop}

\proof As before, the constants $\beta_{\bar{x},i},\beta_{\bar{y},i},
\beta_{\bar{z},i}$ are the sums of infinite series corresponding to all
immersed triangular regions with vertices at $a_i$, $c_i$, and one of
$\bar{x},\bar{y},\bar{z}$. For example the coefficient
$\beta_{\bar{z},i}$ associated to composition
$m_2(\bar{x},c_i)$ is the sum of an infinite series of contributions
associated to triangular regions $T_{\bar{z},i,n}$, $n\in\Z$. The integer
labels can be chosen so that $\partial T_{\bar{z},i,n}-
\partial T_{\bar{z},i,n'}=(n-n')\partial C_i$ and, after deforming to the
configuration in Figure \ref{fig:deltaij}, $T_{\bar{z},i,0}$ becomes
the smallest triangle with vertices $\bar{z},\,a_i,\,c_i$ (i.e., the
triangle which appears with coefficient $-2$ in the 2-chain $C_i$).
Then one easily checks that $T_{\bar{z},i,n}=T_{\bar{z},i,0}+nC_i+\frac{1}{2}
n(3n-1)\Sigma_0$. Therefore, denoting by $\psi_{\bar{z},i}$ the
coefficient associated to $T_{\bar{z},i,0}$, the same argument as in the
previous calculations yields the formula
$$\beta_{\bar{z},i}=\zeta_{i,-}\psi_{\bar{z},i}.$$

Similarly, denote by $T_{\bar{x},i,n}$, $n\in\Z$, the immersed triangles
contributing to $m_2(\bar{x},c_i)$, in such a way that $\partial T_{\bar{x},
i,n}-\partial T_{\bar{x},i,n'}=(n-n')\partial C_i$, and $T_{\bar{x},i,0}$
corresponds to the smallest triangle with vertices $\bar{x},\,a_i,\,c_i$ in
Figure \ref{fig:deltaij} (i.e.\ the triangle which appears with coefficient
$+2$ in the 2-chain $C_i$). Then $T_{\bar{x},i,n}=T_{\bar{x},i,0}+nC_i+
\frac{1}{2}n(3n+1)\Sigma_0$. Therefore, denoting by $\psi_{\bar{x},i}$ the
contribution of $T_{\bar{x},i,0}$, we have
$\beta_{\bar{x},i}=\zeta_{i,+}\psi_{\bar{x},i}.$

Finally, labelling the triangles with vertices $\bar{y},\,a_i,\,c_i$ by
integers in such a way that $T_{\bar{y},i,0}$ and
$T_{\bar{y},i,1}=T_{\bar{y},i,0}+C_i$ correspond to the negative and
positive parts of $C_i$ respectively, it is easy to check that
$T_{\bar{y},i,n}=T_{\bar{y},i,0}+nC_i+\frac{3}{2}n(n-1)\Sigma_0$, so denoting
by $\psi_{\bar{y},i}$ the contribution of $T_{\bar{y},i,0}$ we have
$\beta_{\bar{y},i}=\zeta_{i,0}\psi_{\bar{y},i}.$
It follows that
$$\frac{\beta_{\bar{z},i}^2\,\alpha_{xy}\alpha_{zz}}
{\beta_{\bar{x},i}\beta_{\bar{y},i}\,\alpha_{zy}\alpha_{xz}}
=\frac{\psi_{\bar{z},i}^2\,\psi_{xy}\psi_{zz}}
{\psi_{\bar{x},i}\psi_{\bar{y},i}\,\psi_{zy}\psi_{xz}}\,
\frac{\zeta_{i,-}^2\,\zeta_+\,\zeta_0}
{\zeta_{i,+}\,\zeta_{i,0}\,\zeta_-^2}.$$

Moreover, the 2-chains $2\,T_{\bar{z},i,0}+T_{xy,0}+T_{zz,0}$ and
$T_{\bar{x},i,0}+T_{\bar{y},i,0}+T_{zy,0}+T_{xz,0}$ are equal, which
implies that $\psi_{\bar{z},i}^2\,\psi_{xy}\psi_{zz}=
\psi_{\bar{x},i}\psi_{\bar{y},i}\,\psi_{zy}\psi_{xz}$ and completes
the proof of the first identity.

The arguments are the same for
$$\frac{\beta_{\bar{x},i}^2\,\alpha_{yz}\alpha_{xx}}
{\beta_{\bar{y},i}\beta_{\bar{z},i}\,\alpha_{xz}\alpha_{yx}}
=\frac{\psi_{\bar{x},i}^2\,\psi_{yz}\psi_{xx}}
{\psi_{\bar{y},i}\psi_{\bar{z},i}\,\psi_{xz}\psi_{yx}}\,
\frac{\zeta_{i,+}^2\,\zeta_+\,\zeta_0}
{\zeta_{i,0}\,\zeta_{i,-}\,\zeta_-^2},
$$
observing that
$2\,T_{\bar{x},i,0}+T_{yz,0}+T_{xx,0}-T_{\bar{y},i,0}-T_{\bar{z},i,0}-
T_{xz,0}-T_{yx,0}=C_i$ (for which the corresponding weight is $-q_i$),
and for
$$\frac{\beta_{\bar{y},i}^2\,\alpha_{zx}\alpha_{yy}}
{\beta_{\bar{z},i}\beta_{\bar{x},i}\,\alpha_{yx}\alpha_{zy}}
=\frac{\psi_{\bar{y},i}^2\,\psi_{zx}\psi_{yy}}
{\psi_{\bar{z},i}\psi_{\bar{x},i}\,\psi_{yx}\psi_{zy}}\,
\frac{\zeta_{i,0}^2\,\zeta_+\,\zeta_0}
{\zeta_{i,-}\,\zeta_{i,+}\,\zeta_-^2}
,$$
observing that
$2\,T_{\bar{y},i,0}+T_{zx,0}+T_{yy,0}-T_{\bar{z},i,0}-T_{\bar{x},i,0}-
T_{yx,0}-T_{zy,0}=\Sigma_0-C_i$ (for which the corresponding weight is
$-q_F/q_i$).
\endproof

\begin{cor} The constants
$\beta_{\bar{x},i},\,\beta_{\bar{y},i},\,\beta_{\bar{z},i}$
satisfy the relations:
\ $\displaystyle \frac{\beta_{\bar{z},i}^3}{\beta_{\bar{x},i}^3}
\,\frac{\alpha_{xy}\alpha_{yx}\alpha_{zz}}{\alpha_{yz}\alpha_{zy}\alpha_{xx}}
= -\frac{1}{q_i}\Bigl(\frac{\zeta_{i,-}}{\zeta_{i,+}}\Bigr)^{\!3}$,
\ $\displaystyle \frac{\beta_{\bar{x},i}^3}{\beta_{\bar{y},i}^3}
\,\frac{\alpha_{yz}\alpha_{zy}\alpha_{xx}}{\alpha_{zx}\alpha_{xz}\alpha_{yy}}
= \frac{q_i^2}{q_F}\Bigl(\frac{\zeta_{i,+}}{\zeta_{i,0}}\Bigr)^{\!3}$, \ and
\ $\displaystyle \frac{\beta_{\bar{y},i}^3}{\beta_{\bar{z},i}^3}
\,\frac{\alpha_{zx}\alpha_{xz}\alpha_{yy}}{\alpha_{xy}\alpha_{yx}\alpha_{zz}}
= -\frac{q_F}{q_i}\Bigl(\frac{\zeta_{i,0}}{\zeta_{i,-}}\Bigr)^{\!3}.$
\end{cor}

\begin{prop}\label{prop:m2delta}
For all $0\le i,j<k$ we have the identities
$$
\frac{\beta_{\bar{y},i}\,\beta_{\bar{z},j}}{\beta_{\bar{y},j}\,
\beta_{\bar{z},i}}=\tilde{q}_{i,j}\,\frac{\zeta_{i,0}\,\zeta_{j,-}}
{\zeta_{j,0}\,\zeta_{i,-}}, \ \
\frac{\beta_{\bar{z},i}\,\beta_{\bar{x},j}}{\beta_{\bar{z},j}\,
\beta_{\bar{x},i}}=\tilde{q}_{i,j}\,\frac{\zeta_{i,-}\,\zeta_{j,+}}
{\zeta_{j,-}\,\zeta_{i,+}},\ \ \mathit{and}\ \
\frac{\beta_{\bar{x},i}\,\beta_{\bar{y},j}}{\beta_{\bar{x},j}\,
\beta_{\bar{y},i}}=\tilde{q}_{i,j}^{-2}\,\frac{\zeta_{i,+}\,\zeta_{j,0}}
{\zeta_{j,+}\,\zeta_{i,0}},$$
where $\tilde{q}_{i,j}=\exp(2\pi i\,[B+i\omega]\cdot[\bar\Delta_{i,j}])$, and
$\zeta_{i,+},\,\zeta_{i,-},\,\zeta_{i,0}$ are
as in Definition \ref{def:zetais}.
\end{prop}

\proof We claim that $T_{\bar{y},i,0}+T_{\bar{z},j,0}-T_{\bar{y},j,0}-
T_{\bar{z},i,0}=\Delta_{i,j}$. Indeed, consider first a situation in
which $L_{3+i}$ lies in the position represented in Figure \ref{fig:deltaij},
and $L_{3+j}$ lies close to it, but is slightly shifted towards the
lower-right direction. Then the intersection points $a_j$ and $c_j$ lie
close to $a_i$ and $c_i$, and following the triangular regions through the
small deformation which takes $L_{3+i}$ to $L_{3+j}$, we easily see that
$T_{\bar{z},j,0}$ is obtained by slightly truncating $T_{\bar{z},i,0}$ on
its $L_{3+i}$ side. Similarly, $T_{\bar{y},j,0}$ is obtained by slightly
truncating $T_{\bar{y},i,0}$, and since $\Delta_{i,j}$ is simply the thin
strip in between $L_{3+i}$ and $L_{3+j}$ the claim follows.

The same property remains true if $L_{3+i}$ and $L_{3+j}$ are further apart
from each other. This can be checked explicitly for example in the
configuration of Figure \ref{fig:sigma0}, where $\Delta_{i,j}$ is as
pictured on Figure~\ref{fig:deltaij} (right). (In this configuration the
deformation from $L_{3+i}$ to $L_{3+j}$ passes through $\bar{y}$ and
$\bar{z}$, so the triangles $T_{\bar{z},i,0}$ and $T_{\bar{z},j,0}$ lie on
opposite sides of $\bar{z}$, and similarly for $T_{\bar{y},i,0}$ and
$T_{\bar{y},j,0}$; this latter triangle is now the small region to the
lower-right of $\bar{y}$ on Figure~\ref{fig:sigma0}).

As a consequence, we have the identity
$$\frac{\psi_{\bar{y},i}\psi_{\bar{z},j}}{\psi_{\bar{y},j}
\psi_{\bar{z},i}}=\tilde{q}_{i,j},$$
which implies the first formula in the proposition. The two other formulas
are proved similarly, using the equalities $T_{\bar{z},i,0}+T_{\bar{x},j,0}-
T_{\bar{z},j,0}-T_{\bar{x},i,0}=\Delta_{i,j}$ and
$T_{\bar{x},i,0}+T_{\bar{y},j,0}-
T_{\bar{x},j,0}-T_{\bar{y},i,0}=-2\,\Delta_{i,j}$.
\endproof

\begin{remark}
The various ratios computed in Propositions \ref{prop:m2p2}--\ref{prop:m2delta}
are {\it intrinsic} quantities attached to the symplectic geometry of $W_k$,
i.e.\ they are invariant under Hamiltonian deformations, irrespective of
whether the vanishing cycles are represented by geodesics or not.
Equivalently, they are invariant under {\it rescalings} of the chosen generators of the morphism
spaces in $\FS(W_k,\{\gamma_i\})$. On the other hand, if we allow ourselves
to use the fact that the vanishing cycles are geodesics in a flat torus, we
can also compute some interesting non-intrinsic quantities (i.e., quantities
which depend on a particular choice of scaling of the generators).

For example, the invariance of $L_0,L_1,L_2$ under the translation of the
torus which maps $x_0$ to $y_0$ (and $y_0$ to $z_0$, $z_0$ to $x_0$) implies
that, for suitable choices of the marked points associated to the spin
structures and of the isomorphisms between lines used to calculate boundary
holonomies,
$\alpha_{xy}=\alpha_{yz}=\alpha_{zx}$,
$\alpha_{yx}=\alpha_{xy}=\alpha_{xz}$, and 
$\alpha_{xx}=\alpha_{yy}=\alpha_{zz}$. In fact, going over the calculations
in the proofs of Propositions \ref{prop:m2p2} and \ref{prop:m2p2f}, and
observing that, in terms of areas and boundary holonomies, 
the contributions of $T_{xy,0}-T_{yx,0}$ and $T_{xx,0}-T_{yx,0}$
are equivalent to those of $\frac13 C$ and $\frac13(\Sigma_0-C)$
respectively, one
easily checks that there exists a constant $s\neq 0$ such that
\begin{equation}\label{eq:comp1}
\begin{array}{l}
\alpha_{xy}=\alpha_{yz}=\alpha_{zx}=s\,q_C^{1/3}\,\zeta_+,\\
\alpha_{xx}=\alpha_{yy}=\alpha_{zz}=s\,q_F^{1/3}q_C^{-1/3}\,\zeta_0,\\
\alpha_{yx}=\alpha_{zy}=\alpha_{xz}=-s\,\zeta_-,
\end{array}
\end{equation}
where by definition $q_C^{1/3}=\exp(\frac{2\pi i}{3}
[B+i\omega]\cdot[\bar{C}])$ and $q_F^{1/3}=\exp(\frac{2\pi i}{3}\,[B+i\omega]\cdot[\Sigma_0])$.
Similarly, for suitable choices we have
\begin{equation}\label{eq:comp2}
\beta_{\bar{x},i}=s_i\,q_i^{1/3}\,\zeta_{i,+},\quad
\beta_{\bar{y},i}=s_i\,q_F^{1/3}\,q_i^{-1/3}\,\zeta_{i,0},\quad
\mathrm{and}\ \beta_{\bar{z},i}=-s_i\,\zeta_{i,-},
\end{equation}
where $s_i$ is a non-zero constant and $q_i^{1/3}=\exp(\frac{2\pi i}{3}
[B+i\omega]\cdot[\bar{C}_i])$. 

The formulas (\ref{eq:comp1}) and
(\ref{eq:comp2}) are only valid in the flat case, when the
complexified symplectic form on $\Sigma_0$ is translation-invariant and
the vanishing cycles are geodesics; however, in the general case we can
always modify our choices of generators of the various morphism spaces
by suitable scaling factors (or equivalently, modify the vanishing cycles by
certain Hamiltonian isotopies) in order to make these formulas hold.
It is therefore these simpler formulas that we will use in order to
determine the mirror map in \S \ref{s:mainproof} below.
\end{remark}

\subsection{Simple degenerations}\label{ss:simpdegsympl}

In this section we consider the situation where the symplectic area of
one of the 2-cycles $\bar\Delta_{i,j}$ becomes a multiple of that of the
fiber $\Sigma_0$. The vanishing cycles $L_{3+i}$ and $L_{3+j}$ are then
Hamiltonian isotopic to each other in $\Sigma_0$, and hence cannot be
represented by disjoint geodesics anymore. However we can still represent $L_{3+i}$
by a closed geodesic, and $L_{3+j}$ by a small generic Hamiltonian
perturbation of $L_{3+i}$, intersecting it transversely in two points.
These two intersection points have Maslov indices $0$ and $1$ respectively
(if we choose the same graded lifts as previously), and for this
configuration we have:

\begin{lemma}
If there exist integers $n\in\Z$ and $i<j$ such that
$[\omega]\cdot [\bar\Delta_{i,j}]=n\,[\omega]\cdot [\Sigma_0]$,
then $\Hom(L_{3+i},L_{3+j})$ is graded isomorphic to
$H^*(S^1)\otimes\C$. Moreover, the differential $$m_1:\Hom^0(L_{3+i},
L_{3+j})\to \Hom^1(L_{3+i},L_{3+j})$$ is zero if $[B]\cdot [\bar\Delta_{i,j}]
\in \Z+n\,[B]\cdot[\Sigma_0]$, and an isomorphism otherwise.
\end{lemma}

\proof The only contributions to $m_1$ come from the two disks $D'$ and $D''$
bounded by $L_{3+i}$ and $L_{3+j}$. The 2-chain $D'-D''$ in $\Sigma_0$ has
symplectic area zero, and is in fact given by
$D'-D''=\Delta_{i,j}-n\Sigma_0$. Hence we can compare the coefficients
$\psi'$ and $\psi''$ associated to these two disks by the same argument
as in \S \ref{ss:m2}. Namely, $\psi'$ and $\psi''$ differ by a sign
factor, a holonomy factor, and an area factor.

In this case the sign factor is $-1$ (the sign rule for odd degree
morphisms is slightly more subtle than that for even degree morphisms
\cite{SeBook}; here we can
see directly that the signs for $D'$ and $D''$ have to be different since
the untwisted Floer homology of $L_{3+i}$ and $L_{3+j}$ is non-trivial);
the holonomy factor is the total holonomy along
$\partial(D'-D'')=L_{3+j}-L_{3+i}$, i.e.\ $\exp(2\pi i\int_{D_{3+i}-D_{3+j}} B)$;
and the area factor is $\exp(2\pi i\int_{D'-D''} B+i\omega)$. It follows that
$$\psi'=-\exp\bigl(2\pi i\,[B+i\omega]\cdot([\bar\Delta_{i,j}]-n[\Sigma_0])
\bigr)\,\psi'',$$
since $D'-D''+D_{3+i}-D_{3+j}=\bar{\Delta}_{i,j}-n\Sigma_0$.
Since $m_1$ is determined by the sum $\psi'+\psi''$, we conclude that
$m_1=0$ if and only if $[B+i\omega]\cdot([\bar\Delta_{i,j}]-n[\Sigma_0])$
is an integer.
\endproof

In other words, if $[B+i\omega]\cdot[\bar{\Delta}_{i,j}]\in \Z\oplus
([B+i\omega]\cdot[\Sigma_0])\,\Z$, then $(L_{3+i},\nabla_{3+i})$ and
$(L_{3+j},\nabla_{3+j})$ are essentially identical, and we have a
non-cancelling pair of extra morphisms of degrees $0$ and $1$ from $L_{3+i}$
to $L_{3+j}$; this mirrors the situation in which $\CP^2$ is blown up twice
at infinitely close points, in which case there is a rational $-2$-curve and
the derived category of coherent sheaves is richer than in the
generic case. In all other situations the intersection points between
$L_{3+i}$ and $L_{3+j}$, if any, are killed by the twisted Floer
differential (even when $L_{3+i}$ and $L_{3+j}$ are Hamiltonian isotopic).

\begin{remark} 
It is important to note that, due to 
the presence of immersed convex polygonal regions with two edges on $L_0\cup
L_1\cup L_2$ and two edges on $L_{3+i}\cup L_{3+j}$ (with a corner at the
intersection point of Maslov index $1$), we have to consider
not only the Floer differential $m_1$, but also the
higher-order composition $m_3$. For example, when $L_{3+i}$ and $L_{3+j}$
are Hamiltonian isotopic the
composition $$m_3:\Hom(L_0,L_2)\otimes \Hom(L_2,L_{3+i})\otimes
\Hom^1(L_{3+i},L_{3+j})\lto \Hom(L_0,L_{3+j})$$ is in general non-zero
(and similarly with $L_1$ instead of $L_0$ or $L_2$).
\end{remark}

As in \S \ref{ss:simpdeg}, it is possible to describe things in a simpler
and more
unified manner by considering a suitable mutation of the exceptional
collection $(L_0,\dots,L_{k+2})$. Assume for simplicity that the two
vanishing cycles which may coincide are $L_3$ and $L_4$, while the others
are represented by distinct geodesics. Then we can modify
the system of arcs $\{\gamma_i\}$ considered so far to a new ordered 
system of arcs  $\{\gamma'_i\}$ such that $\gamma'_i=\gamma_i$ for
$i\not\in\{2,3\}$, $\gamma'_3=\gamma_2$, and $\gamma'_2$ connects the
origin to $\lambda_3\approx \infty$ along the negative real axis. This
gives rise to a new category $\FS(W_k,\{\gamma'_i\})$, in which all objects
but one can be identified with the objects $L_i,\ i\neq 3$ of $\FS(W_k,\{\gamma_i\})$; thus,
we denote by $L_0,L_1,L',L_2,L_4,\dots,L_{k+2}$ the objects of
$\FS(W_k,\{\gamma'_i\})$. The morphisms and compositions not involving $L'$
are as in $\FS(W_k,\{\gamma_i\})$. 

The new vanishing cycle $L'$ is
Hamiltonian isotopic to the image of $L_3$ under the positive Dehn twist
along $L_2$. In particular, with the notations of Lemma \ref{l:vch1}, and
for a suitable choice of orientation, its
homology class is $[L']=[L_2]-[L_3]=b$. Choosing a geodesic representative, we
have $|L_0\cap L'|=2$, $|L_1\cap L'|=1$, $|L'\cap L_2|=1$, and $|L'\cap
L_{3+i}|=1$ for $i\ge 1$, and all morphisms in $\FS(W_k,\{\gamma'_i\})$ have
degree 0.

Because $L'$ is Hamiltonian isotopic to the image of $L_3$ under the Dehn
twist along $L_2$, the fiber $\Sigma_0$ contains a 2-chain $\Delta'$ with 
$\partial \Delta'=L'+L_4-L_2$ and such that $\int_{\Delta'}\omega=
\int_{\Delta_{3,4}}\omega$. Capping off $\Delta'$ with the appropriate
Lefschetz thimbles, we obtain a 2-cycle $\bar\Delta'$ in $M_k$, with
$[\bar\Delta']=[\bar\Delta_{3,4}]$ in $H_2(M_k,\Z)$.
The composition $$\Hom(L',L_2)\otimes \Hom(L_2,L_4)\lto \Hom(L',L_4)$$
corresponds to an infinite series of triangular immersed regions in
$\Sigma_0$, of which in general two are embedded. The case where the
symplectic area of $\Delta'$ is a multiple of that of the fiber corresponds
precisely to the situation where the two embedded triangular regions have
equal symplectic areas. In general, the immersed triangles contributing to
the composition can be labelled $T'_n$, $n\in\Z$, in such a way that
$T'_n=T'_0+n\Delta'+\frac12 n(n-1) \Sigma_0$. Arguing as before, one
easily shows that the composition is given by the contribution of $T'_0$
multiplied by the factor
$$\sum_{n\in\Z} (-1)^n\,{q'}^{\,n}\,q_F^{n(n-1)/2},\quad
\mathrm{where}\ q'=\exp(2\pi i[B+i\omega]\cdot[\bar\Delta'])=\tilde{q}_{3,4}.$$
This multiplicative factor vanishes if and only if $q'=q_F^k$ for some $k\in
\Z$ (an easy way to see this is to view this factor as a theta
function, see below), i.e.\ iff
$[B+i\omega]\cdot[\bar\Delta']\in
\Z\oplus ([B+i\omega]\cdot[\Sigma_0])\Z$. Hence, as in \S \ref{ss:simpdeg}
the mutation makes it possible to avoid dealing with a
non-trivial differential, and provides an alternative description in which
the simple degeneration corresponds to one of the composition maps becoming
identically zero.

\subsection{Modular invariance and theta functions} \label{ss:theta}

In this section we study the modularity properties of the category
$\FS(W_k,\{\gamma_i\})$ with respect to some of the parameters governing
deformations of the complexified symplectic structure, and the relation with
theta functions.

\begin{prop}\label{prop:periodic}
Consider two complexified symplectic forms $\kappa=B+i\omega$ and
$\kappa'=B'+i\omega'$ on $M_k$, such that
$[\kappa']\cdot[\Sigma_0]=[\kappa]\cdot[\Sigma_0]$ and $[\kappa']-[\kappa]
\in H^2(M_k,\Z)\oplus (\kappa\cdot[\Sigma_0])\,H^2(M_k,\Z)$.
Then the categories $\FS(W_k,\kappa,\{\gamma_i\})$ and $\FS(W_k,\kappa',
\{\gamma_i\})$ are equivalent.
\end{prop}

\proof
First consider the situation where $\omega'=\omega$, and $B'=B+d\chi$ for
some 1-form $\chi$. Then the vanishing cycles $L_i$ remain the same, but
the associated flat connections differ, and we can e.g.\ take
$\nabla'_i=\nabla_i-2\pi i\chi$. Then the contribution of a pseudo-holomorphic
map $u:(D^2,\partial D^2)\to(\Sigma_0,\cup L_i)$ is actually the same in
both cases, since the holonomy term changes by $\exp(-2\pi i\int_{u(\partial
D^2)}\chi)$, while the weight factor changes by
$\exp(2\pi i\int_{D^2}u^*d\chi)=\exp(2\pi i\int_{u(\partial D^2)}\chi)$.
So, in the more general situation where $[B'-B]\in H^2(M_k,\Z)$ and
$[B'-B]\cdot[\Sigma_0]=0$ (still assuming $\omega'=\omega$), after
modifying $B$ by an exact term we can assume that $B$ and $B'$ coincide
over $\Sigma_0$, and that the integral of $B'-B$ over each thimble $D_i$
is a multiple of $2\pi$. In this situation the vanishing cycles $L_i$
are the same, and the associated flat connections are gauge equivalent
(since their holonomies differ by multiples of $2\pi$), so the corresponding
twisted Floer theories are identical.

Next, consider the situation where $[B+i\omega]$ changes by an integer
multiple of $[B+i\omega]\cdot[\Sigma_0]$. After adding an exact term to
$\kappa=B+i\omega$ (which does not affect the category of vanishing cycles by Lemma
\ref{l:exactdeform} and by the above remark), we can assume that $\kappa$
and $\kappa'$ coincide over $\Sigma_0$, and that the relative cohomology class
of $\kappa'-\kappa$ is an element of
$(\kappa\cdot[\Sigma_0])H^2(M_k,\Sigma_0;\Z)$.

Let $D_i$ and $D'_i$ be the thimbles associated to the arc $\gamma_i$ and
to the symplectic forms $\omega$ and $\omega'$ respectively. The integrality
assumption on $\kappa'-\kappa$ implies that there exists an integer
$n_i\in\Z$ such that $\int_{D_i} \kappa'=n_i [\kappa]\cdot[\Sigma_0]+
\int_{D_i} \kappa$. Since $D_i$ and $D'_i$ can be deformed
continuously into each other (by deforming the horizontal distribution),
there exists a 2-chain $K_i$ in $\Sigma_0$ such that $[D_i+K_i-D'_i]=0$ in
$H_2(M_k)$. Then $\int_{K_i}\omega=\int_{K_i}\omega'=
-\int_{D_i}\omega'=-n_i[\omega]\cdot[\Sigma_0]$.
Since the symplectic area of the 2-chain $K_i\subset\Sigma_0$ is an integer
multiple of that of the fiber,
the two vanishing cycles $L'_i=\partial D'_i$ and $L_i=\partial D_i$ are
mutually Hamiltonian isotopic in $\Sigma_0$, and hence
we can assume that $L'_i=L_i$. Moreover,
in $H_2(M_k,L_i)$ we have $[D'_i]=[D_i]-n_i[\Sigma_0]$.
Therefore, $\int_{D'_i}B'=\int_{D_i}B'-n_i\int_{\Sigma_0}B'=(\int_{D_i}B+
n_i[B]\cdot[\Sigma_0])-n_i[B]\cdot[\Sigma_0]=\int_{D_i}B$. So the flat
connections $\nabla_i$ and $\nabla'_i$ have the same holonomy, which implies
that $(L_i,\nabla_i)$ and $(L'_i,\nabla'_i)$ behave identically for twisted
Floer theory.
\endproof

This property explains the invariance of the structure coefficients
($\alpha_{xy}$, etc.) under certain changes of variables. More precisely,
one easily checks that 
$\zeta_+(q_C q_F^3,q_F)=-q_C^{-1}q_F^{-2}
\zeta_+(q_C,q_F)$, $\zeta_-(q_C q_F^3,q_F)=-q_C^{-1}q_F^{-1}\zeta_-(q_C,q_F)$,
and $\zeta_0(q_C q_F^3,q_F)=-q_C^{-1}\zeta_0(q_C,q_F)$. This implies that
the quantities considered in
Propositions \ref{prop:m2p2} and \ref{prop:m2p2f} are invariant under the
change of variables $(q_C,q_F)\mapsto (q_C q_F^3,q_F)$; a closer examination
shows that the individual constants $\alpha_{xy}$, etc.\ are also invariant
under this change of variables.

On the other hand, one easily checks that $\zeta_+(q_Cq_F,q_F)=-q_C^{-1}\zeta_0(q_C,q_F)$,
$\zeta_0(q_Cq_F,q_F)=\zeta_-(q_C,q_F)$, and $\zeta_-(q_Cq_F,q_F)=
\zeta_+(q_C,q_F)$, which may seem surprising at first. The reason is
that this change of variables corresponds to a
non-Hamiltonian deformation of e.g.\ $L_1$ which sweeps exactly
once through the entire fiber $\Sigma_0$. This deformation preserves the
intersection points, but induces a non-trivial permutation of their labels:
namely, $x_0,\,y_0,\,z_0$ become $y_0,\,z_0,\,x_0$ respectively, and
$x_1,\,y_1,\,z_1$ become $z_1,\,x_1,\,y_1$ respectively. Thus, for example,
$\alpha_{xy}(q_C,q_F)=\alpha_{yx}(q_Cq_F,q_F)=\alpha_{zz}(q_Cq_F^2,q_F)$
(and similarly for the other coefficients).

\medskip

Another way to understand these invariance properties is to relate the
functions $\zeta_+$, $\zeta_-$, and $\zeta_0$ to theta functions.
Recall that the ordinary theta function is an analytic function defined by
$$
\theta(z,\tau)=\sum_{n\in\Z}\exp(\pi i n^2\tau+2\pi i nz),
$$
where $z\in\C$ and $\tau\in \cH$ (here $\cH$ is the upper half-plane
$\{\mathrm{Im}\,\tau>0\}$).
This function is quasiperiodic with respect to
the lattice $\Lambda_{\tau}\subset \C$ generated by $1$ and $\tau$, and its
behavior under translation by an element of the lattice
is given by the formula
$$
\theta(z+u\tau+v, \tau)=\exp(-\pi i u^2\tau-2\pi i u z)\theta(z, \tau).
$$
The zeros of the theta function are the infinite set
$\left\{z=(n+\frac12)+(m+\frac12)\tau\; |\; n,m\in \Z\right\}.$

Here we consider theta functions with rational characteristics
$a,b\in \mathbb{Q}$, defined by
$$
\theta_{a,b}(z,\tau)=\sum_{n\in\Z}\exp(\pi i (n+a)^2\tau+2\pi
i(n+a)(z+b)).
$$
Let us introduce   new variables $q=\exp(\pi i\tau)$ and $w=\exp(\pi iz).$
Now the following three $\theta$-functions play a very important role in our
considerations:
\begin{eqnarray*}
\theta_{\frac{1}{2},\frac{1}{2}}(3z,3\tau)&=&\exp(\tfrac{i\pi}{2})\,q^{3/4}\,
\sum_{n\in \Z}(-1)^n w^{6n+3}q^{3n^2+3n},\\
\theta_{\frac{1}{6},\frac{1}{2}}(3z,3\tau)&=&\exp(\tfrac{i\pi}{6})\,q^{1/12}\,
\sum_{n\in \Z} (-1)^nw^{6n+1}q^{3n^2+n},\\
\theta_{\frac{5}{6},\frac{1}{2}}(3z,3\tau)&=&\exp(-\tfrac{i\pi}{6})\,q^{1/12}\,
\sum_{n\in \Z} (-1)^nw^{6n-1}q^{3n^2-n}.
\end{eqnarray*}
The zero set of the function $\theta_{\frac{1}{2},\frac{1}{2}}(3z,3\tau)$
is $\left\{\frac{n}{3}+m\tau \;|\; n,m\in \Z\right\}$, while the zero sets
of the functions
$\theta_{\frac{1}{6},\frac{1}{2}}(3z,3\tau)$ and $\theta_{\frac{5}{6},\frac{1}{2}}(3z,3\tau)$ are
$$\left\{\tfrac{n}{3} + (m+\tfrac{1}{3})\tau \;|\; n,m\in \Z\right\}\quad \text{and}
\quad\left\{\tfrac{n}{3} + (m-\tfrac{1}{3})\tau\;|\; n,m\in \Z\right\}
$$ respectively.
These three theta functions can be viewed as holomorphic sections 
of a line bundle of degree 3
on the elliptic curve $E=\C/\Lambda_{\tau}$; considering the zero sets,
we see that this line bundle is $\mathbb{L}=\O_E(3\cdot (0))$.
These three sections of $\mathbb{L}$ determine an embedding of the elliptic
curve $E=\C/\Lambda_{\tau}$ into the projective plane, given by
$$
z\mapsto (\theta_{\frac{1}{2},\frac{1}{2}}(3z,3\tau):
\theta_{\frac{1}{6},\frac{1}{2}}(3z,3\tau):
\theta_{\frac{5}{6},\frac{1}{2}}(3z,3\tau)).
$$
Observe that the two functions $$\theta_{\frac{1}{2},\frac{1}{2}}(3z,3\tau)\,
\theta_{\frac{1}{6},\frac{1}{2}}(3z,3\tau)
\theta_{\frac{5}{6},\frac{1}{2}}(3z,3\tau)\quad\mathrm{and}\quad
\theta_{\frac{1}{2},\frac{1}{2}}(3z,3\tau)^3+
\theta_{\frac{1}{6},\frac{1}{2}}(3z,3\tau)^3+
\theta_{\frac{5}{6},\frac{1}{2}}(3z,3\tau)^3$$
coincide up to a constant multiplicative factor, since they both
correspond to holomorphic sections of the line bundle $\mathbb{L}^{\otimes
3}$ over $E$, and an easy calculation shows that they
have the same zero set $\{\frac{n}{3}+\frac{m}{3}\tau\,|\,n,m\in\Z\}$.
Therefore, the image of the above embedding of $E$ into $\PP^2$
is the cubic given by the equation
$$
(A^3+B^3+C^3)XYZ-ABC(X^3+Y^3+Z^3)=0,
$$
where $(A,B,C)$ are the values of the three theta functions at any given
point of $\C/\Lambda_{\tau}$ (not in $\frac13 \Lambda_\tau$).

Consider the function
$$
\left(\frac{\theta_{\frac{1}{6},\frac{1}{2}}(3z,3\tau)}{\theta_{\frac{5}{6},\frac{1}{2}}(3z,3\tau)}\right)^{\!3}=
-\left(\frac{\sum_{n\in\Z}\, (-1)^n \,w^{6n+1}\,
q^{n(3n+1)}}
{\sum_{n\in\Z}\, (-1)^n\, w^{6n-1}\, q^{n(3n-1)}}\right)^{\!3}.
$$
Substituting $q^2=q_F$ and $w^6=q_C$, one easily checks that this
coincides with the expression which appears in Proposition \ref{prop:m2p2},
$$
\frac{\alpha_{xy}\alpha_{yz}\alpha_{zx}}{\alpha_{yx}\alpha_{zy}
\alpha_{xz}}=-q_C\left(\frac{\sum_{n\in\Z}\, (-1)^n \,q_C^{n\vphantom{/}}\,
q_F^{n(3n+1)/2}
} {\sum_{n\in\Z}\, (-1)^n\, q_C^{n\vphantom{/}}\, q_F^{n(3n-1)/2}
}\right)^{\!3}.
$$
Similarly,
$$
\left(\frac{\theta_{\frac{1}{2},\frac{1}{2}}(3z,3\tau)}{\theta_{\frac{5}{6},\frac{1}{2}}(3z,3\tau)}\right)^{\!3}=
q^2\left(\frac{\sum_{n\in \Z}(-1)^n w^{6n+3}q^{3n^2+3n}}{\sum_{n\in \Z}
(-1)^nw^{6n-1}q^{3n^2-n}}\right)^{\!3}=
-q^2\left(\frac{\sum_{n\in \Z}(-1)^n w^{6n-3}q^{3n^2-3n}}{\sum_{n\in \Z}
(-1)^nw^{6n-1}q^{3n^2-n}}\right)^{\!3}.
$$
After the same substitution $q^2=q_F$ and $w^6=q_C$, this coincides
with the expression given in Proposition \ref{prop:m2p2f},
$$\displaystyle
\frac{\alpha_{xx}\alpha_{yy}\alpha_{zz}}{\alpha_{yx}\alpha_{zy}
\alpha_{xz}}=-\frac{q_F}{q_C}
\left(\frac{\sum_{n\in\Z}\, (-1)^n \,q_C^{n\vphantom{/}}\,
q_F^{3n(n-1)/2}
} {\sum_{n\in\Z}\, (-1)^n\, q_C^{n\vphantom{/}}\, q_F^{n(3n-1)/2}
}\right)^{\!3}.$$
Similarly, in the case where (\ref{eq:comp1}) holds, one easily checks that
\begin{equation}\label{eq:alphatheta}\begin{array}{l}
\alpha_{xy}=\alpha_{yz}=\alpha_{zx}=\tilde{s}\,e^{-2i\pi/3}\,
\theta_{\frac16,\frac12}(3z_0,3\tau),\\
\alpha_{xx}=\alpha_{yy}=\alpha_{zz}=\tilde{s}\,
\theta_{\frac12,\frac12}(3z_0,3\tau),\\
\alpha_{yx}=\alpha_{zy}=\alpha_{xz}=\tilde{s}\,e^{2i\pi/3}\,
\theta_{\frac56,\frac12}(3z_0,3\tau),
\end{array}
\end{equation}
where $\tau=[B+i\omega]\cdot [\Sigma_0]$,
$z_0=\frac{1}{3}[B+i\omega]\cdot[\bar{C}]$,
and $\tilde{s}=e^{i\pi/2}\,q_F^{-1/24}\,q_C^{1/6}\,s\neq 0$.
Similar interpretations can be made for the quantities considered in
Propositions \ref{prop:m2beta}--\ref{prop:m2delta} and in (\ref{eq:comp2}).

\section{Proof of the main theorems}\label{s:mainproof}

The derived categories considered in \S \ref{s:dbcoh} depend on 
an elliptic curve $E$,
two degree 3 line bundles $\L_1,\L_2$ over $E$, and $k$ points
$p_1,\dots,p_k$ on $E$. Meanwhile, the
categories considered in \S \ref{s:lagvc} depend on a cohomology class
$[B+i\omega]\in H^2(M_k,\C)$. We now show how to relate these
two sets of parameters.

Fix the cohomology class $[B+i\omega]\in H^2(M_k,\C)$, and consider
the category $\db{\FS(W_k)}$ studied in \S \ref{s:lagvc}. With the notations
of \S \ref{ss:h2mk}, assume that $[\omega]\cdot [\bar\Delta_{i,j}]$ is not
an integer multiple of $[\omega]\cdot[\Sigma_0]$ for any
$i,j\in\{0,\dots,k-1\}$. Then $\db{\FS(W_k)}$ admits a full strong
exceptional collection $(L_0,\dots,L_{k+2})$, whose properties have been 
studied in \S \ref{s:lagvc}. In particular, the objects and morphisms in
this exceptional collection are the same as for the exceptional collection
$\sigma=(\O_{X_K}, \pi^*\T_{\PP^2}(-1), \pi^*\O_{\PP^2}(1), \O_{l_1},\dots,
\O_{l_k})$ considered in \S \ref{s:dbcoh} for the derived category of
coherent sheaves on a (possibly noncommutative) Del Pezzo surface.
Hence, our goal is now to compare the composition laws and show that, 
for a suitable choice of the parameters $(E,\mathcal{L}_1,\mathcal{L}_2,K)$,
the algebra of homomorphisms of the exceptional collection
$(L_0,\dots,L_{k+2})$ is isomorphic to the algebra $B_{K,\mu}$ considered
in \S \ref{s:dbcoh}. More precisely, we claim:

\begin{prop}\label{prop:mirrormap}
Let $E$ be the elliptic curve $\C/\Lambda_\tau$, where
$\tau=[B+i\omega]\cdot[\Sigma_0]$, realized as a plane cubic via
the embedding $j:E\to\PP^2$
given by $z\mapsto (\vartheta_+(z):\vartheta_0(z):\vartheta_-(z))$,
where $$\vartheta_+(z)=e^{-2i\pi/3}\theta_{\frac16,\frac12}(3z,3\tau),\quad
\vartheta_0(z)=\theta_{\frac12,\frac12}(3z,3\tau), \quad \text{and} \ 
\vartheta_-(z)=e^{2i\pi/3}\theta_{\frac56,\frac12}(3z,3\tau).$$ Let $z_0=\frac13
[B+i\omega]\cdot[\bar C]$, and for $i\in\{0,\dots,k-1\}$ let $p_i=
\frac13[B+i\omega]\cdot[\bar C_i]$. Finally, let $\L_1=\O_E(3\cdot(-z_0))$
and $\L_2=\O_E(3\cdot(0))$.
Then the algebra of homomorphisms of the exceptional collection
$(L_0,\dots,L_{k+2})$ is isomorphic to $B_{K,\mu}$, where $\mu$ is
determined by $(E,\L_1,\L_2)$ via Construction \ref{constr:mu}
and $K=\{j(z_0+p_0),\dots,j(z_0+p_{k-1})\}$.
\end{prop}

\proof After a suitable rescaling of the chosen
bases of the morphism spaces (or just by deforming to the situation where
the fiber is flat and the vanishing cycles are geodesics), we can assume
that the compositions of morphisms between the objects $L_0,\dots,L_{k+2}$
are given by the formulas (\ref{eq:comp1}) and (\ref{eq:comp2}).
We identify the vector spaces $U=\Hom(L_0,L_1)$, $V=\Hom(L_1,L_2)$,
and $W=\Hom(L_0,L_2)$ with $\C^3$ by considering the bases
$(x_0,y_0,z_0)$, $(x_1,y_1,z_1)$, and $(\bar{x},\bar{y},\bar{z})$.
The composition tensor $\mu:V\otimes U\to W$ is determined by the
three constants $a=\alpha_{xy}=\alpha_{yz}=\alpha_{zx}$, $b=\alpha_{xx}=
\alpha_{yy}=\alpha_{zz}$, and $c=\alpha_{yx}=\alpha_{zy}=\alpha_{xz}$.
In particular, given an element $v=(X,Y,Z)\in V$, 
the composition map $\mu_v=\mu(v,\cdot):U\to W$ is given by the
matrix \begin{equation}\label{eq:matthetav}\left(\begin{array}{ccc}
\alpha_{xx}X& \alpha_{yz}Z & \alpha_{zy}Y\\
\alpha_{xz}Z& \alpha_{yy}Y & \alpha_{zx}X\\
\alpha_{xy}Y& \alpha_{yx}X & \alpha_{zz}Z
\end{array}\right)=\left(\begin{array}{ccc}
bX& aZ & cY\\
cZ& bY & aX\\
aY& cX & bZ
\end{array}\right)
\end{equation}
which has rank 2 precisely when
\begin{equation}\label{eq:detthetav}
\det(\mu_v)=(a^3+b^3+c^3)\,XYZ-abc(X^3+Y^3+Z^3)=0.
\end{equation}
By (\ref{eq:alphatheta}), the constants $a,b,c$ are (up to a non-zero
constant factor) the values of the theta functions
$\vartheta_+,\vartheta_0,\,\vartheta_-$ at the point $z_0$. Therefore, by the discussion
in \S \ref{ss:theta}, there are two possibilities: 

\begin{enumerate}
\item if $z_0\in\frac13 \Lambda_\tau$, then $abc=0$ and $\mu_v$ always
has rank 2; as explained in \S \ref{ss:noncomm} this corresponds to a
commutative situation;
\item if $z_0\not\in \frac13 \Lambda_\tau$, then (\ref{eq:detthetav})
defines a cubic $\Gamma_V\subset \PP(V)=\PP^2$, and this cubic is precisely the 
image of the embedding $j$.
\end{enumerate}
The same situation holds for $\mu_u$; interestingly, under the chosen
identifications of $\PP(U)$ and $\PP(V)$ with $\PP^2$, the two subschemes
$\Gamma_U\subset \PP(U)$ and $\Gamma_V\subset \PP(V)$ determined by the
equations $\det(\mu_u)=0$ and $\det(\mu_v)=0$ coincide exactly.
However, with this description, the isomorphism $\sigma:\Gamma_V\to\Gamma_U$
which takes $v$ to the point of $\Gamma_U$ corresponding to
$\Ker\,\mu_v$ is {\it not} the identity map. Here the reader
is referred to the discussion on pp.\ 37--38 of \cite{ATV}, which we follow
loosely.

Given a point $v=(X:Y:Z)\in\Gamma_V$, the kernel of $\mu_v$ can
be obtained as the cross-product of any
two of the rows of the matrix (\ref{eq:matthetav}). Taking e.g.\ the
first two rows, we obtain that the corresponding point of $\Gamma_U$ is
\begin{equation}\label{eq:sigma}
\sigma(X:Y:Z)=(a^2 XZ-bc Y^2: c^2 YZ-ab X^2: b^2 XY-ac Z^2).
\end{equation}
Observe that $j$ maps the origin to $(1:0:-1)\in \Gamma_V$, and that the
corresponding point in $\Gamma_U$ is
$\sigma(1:0:-1)=(a:b:c)=j(z_0)$. Hence, considering only the situation where
$\Gamma_U\simeq \Gamma_V\simeq E$, and identifying $E$ with $\Gamma_V$ by
means of the embedding $j$, the identification of $E$ with $\Gamma_U$ is
given by the embedding $\sigma\circ j$, which is the composition of $j$
with the translation by $z_0$. Therefore, the line bundles on $E$
induced by the two inclusions of $E$ into $\PP(U)$ and
$\PP(V)$ are respectively $(\sigma\circ j)^*\O_{\PP^2}(1)=\O_E(3\cdot
(-z_0))=\L_1$ and $j^*\O_{\PP^2}(1)=\O_E(3\cdot (0))=\L_2$. It then follows
from the discussion in \S \ref{ss:noncomm} that the composition tensor
$\mu$ corresponds to the data $(E,\L_1,\L_2)$. This remains true even when
$z_0\in \frac13 \Lambda_\tau$, since in that case we have $\L_1\simeq \L_2$
and the composition tensor associated to the triple $(E,\L_1,\L_2)$ is that
of the usual projective plane (see Remark \ref{rmk:comm}).

Next we consider the composition $\Hom(L_2,L_{3+i})\otimes W\lto \Hom(L_0,
L_{3+i})$. Choosing generators of the lines $\Hom(L_2,L_{3+i})$ and
$\Hom(L_0,L_{3+i})$ we can view this map as a linear form on $W$.
In the given basis of $W$, this
linear form is given by $\bigl(\beta_{\bar{x},i}\,,\,
\beta_{\bar{y},i}\,,\,\beta_{\bar{z},i}\bigr)$,
which by (\ref{eq:comp2}) coincides up to a non-zero constant factor with
$$\bigl(\vartheta_+(p_i)\,,\,\vartheta_0(p_i)\,,\,\vartheta_-(p_i)\bigr).$$
On the other hand we know from \S \ref{ss:noncomm} that the kernel of this
linear form should be exactly $\Im \mu_{v_i}$, where $v_i\in\Gamma_V$ is
the point being blown up.

For any $v=(X:Y:Z)\in\Gamma_V$, the projection $W\to W/\Im \mu_v$ is
a linear form given up to a scaling factor by the dot product of any two
columns of the matrix (\ref{eq:matthetav}). Taking e.g.\ the first two
columns, we obtain that the expression of this linear form 
relatively to our chosen basis of $W$ is 
$$\bigl(c^2 XZ-abY^2\,,\,a^2 YZ - bc X^2 \,,\,b^2 XY - ac Z^2\bigr).$$
Interestingly, if we assume that
$(X:Y:Z)=\sigma(\tilde{X}:\tilde{Y}:\tilde{Z})$, where $\sigma$ is the 
transformation given by (\ref{eq:sigma}), then this expression simplifies
to a scalar multiple of $(\tilde{X}\,,\,\tilde{Y}\,,\,\tilde{Z})$.
Hence, we conclude that $v_i=\sigma(j(p_i))=j(z_0+p_i)$.
\endproof

\begin{remark} \label{rmk:invt}
At this point the reader may legitimately be concerned that, since the
homology classes $[\bar{C}]$ and $[\bar{C}_i]$ are canonically defined only
up to a multiple of $[\Sigma_0]$, and since $[B]$ is only defined up to an
element of $H^2(M_k,\Z)$, the points $z_0$ and $p_i$ of $E$ are
canonically determined only up to translations by elements of $\frac13
\Lambda_\tau$. However, the line bundle $\L_1=\O_E(3\cdot (-z_0))$ is not
affected by this ambiguity in the determination of $z_0$, and neither are
the relative positions of the points $p_i$, since the quantity $p_j-p_i=
[B+i\omega]\cdot [\bar\Delta_{i,j}]$ is well-defined up to an element of
$\Lambda_\tau$. Moreover, a simultaneous translation of all the blown up
points by an element of $\frac13 \Lambda_\tau$ amounts to an automorphism
of the triple $(E,\L_1,\L_2)$, which does not actually affect the category.
(From the point of view of the embedding $j$, this automorphism simply
permutes the homogeneous coordinates $X,Y,Z$ and multiplies them by cubic
roots of unity; this is consistent with the observation made after the proof
of Proposition \ref{prop:periodic}).
\end{remark}

Theorems \ref{th:main1} and \ref{th:main3} now follow directly from the
discussion. Namely, in the case of a blowup of $\CP^2$ at a set
$K=\{p_0,\dots,p_{k-1}\}$ of $k$
distinct points (Theorem \ref{th:main1}), we consider a cubic curve
$E\subset\CP^2$ which contains all the points of $K$, and view it as
an elliptic curve $\C/\Lambda_\tau$ for some $\tau\in\C$ with $\Im\tau>0$.
This allows us to view the points $p_i$ as elements of $\C/\Lambda_\tau$
(well-defined up to a simultaneous translation of all $p_i$ by an element
of $\frac13 \Lambda_\tau$, since the origin can be chosen at any of the
flexes of $E$; however by Remark \ref{rmk:invt} this does not matter for our
construction). Then we equip $M_k$ with a complexified symplectic
structure such that \hbox{$[B+i\omega]\cdot [\Sigma_0]=\tau$}, 
\hbox{$[B+i\omega]\cdot [\bar{C}]=0$}, and 
\hbox{$[B+i\omega]\cdot[\bar{C}_i]=3 p_i$}. The existence of such
a $B+i\omega$ follows from a standard result about symplectic
structures on Lefschetz fibrations:

\begin{prop}[Gompf]
Given any cohomology class $[\zeta]\in H^2(M_k,\R)$ such that
$[\zeta]\cdot [\Sigma_0]>0$, the manifold $M_k$ admits a symplectic
structure in the cohomology class $[\zeta]$, for which the fibers of
$W_k$ are symplectic submanifolds.
\end{prop}

\proof
The map $W_k:M_k\to\C$ is a Lefschetz fibration, and the argument given
in the proof of \cite[Theorem 10.2.18]{GS} can be adapted in a
straightforward manner to this situation, even though the base of the fibration is
not compact. (Alternatively, one can also work with the compactified
fibration $\ol{W_k}:\ol{M}\to\CP^1$).
The symplectic form $\omega$ constructed by this argument lies in the
cohomology class $t[\zeta]+W_k^*([\mathrm{vol}_{\C}])$ for some constant $t>0$;
since the area form on $\C$ is exact, we have $[\omega]=t[\zeta]$, and
scaling $\omega$ by a constant factor we obtain the desired result.
\endproof

By Proposition
\ref{prop:mirrormap} the algebra of homomorphisms of the exceptional
collection $(L_0,\dots,L_{k+2})$ is then isomorphic to $B_K$, which implies
that $\db{\FS(W_k)}\cong \db{\mod B_K}\cong \db{\coh(X_K)}$. 

In the case of
a noncommutative blowup of $\PP^2$ (Theorem \ref{th:main3}), consider the
triple $(E,\L_1,\L_2)$ associated to the underlying noncommutative $\PP^2$,
and view again $E$ as a quotient $\C/\Lambda_\tau$. Choose $z_0$
(well-defined up to an element of $\frac13 \Lambda_\tau$) such that 
$\L_2\otimes \L_1^{-1}\simeq \O_E(3\cdot(z_0)-3\cdot(0))\in \Pic^0(E)$.
As explained in \S \ref{ss:noncomm}, the blown up points must all
lie in $\Gamma_V\subset \PP(V)$, and under the identification
$\Gamma_V\simeq E$ they can be viewed as elements $p_i\in
\C/\Lambda_\tau$. Equip $M_k$ with a complexified symplectic 
structure such that $[B+i\omega]\cdot [\Sigma_0]=\tau$, $[B+i\omega]\cdot
[\bar{C}]=3z_0$, and $[B+i\omega]\cdot[\bar{C}_i]=3 (p_i-z_0)$.
By Proposition \ref{prop:mirrormap} the algebra of homomorphisms of the
exceptional collection $(L_0,\dots,L_{k+2})$ is then isomorphic to
$B_{K,\mu}$, which yields the desired equivalence of categories.

Theorem \ref{th:main2} is proved similarly, working with the mutated
exceptional collections $\tau'$ (introduced in \S \ref{ss:simpdeg}) and
$(L_0,L_1,L',L_2,L_4,\dots,L_{k+2})$ (introduced in \S
\ref{ss:simpdegsympl}). The details are left to the reader.

\begin{remark}
The construction carried out for Theorem \ref{th:main1} also applies to
some limit situations in which $X_K$ is actually not a Del Pezzo surface.
For example, the argument applies equally well to the situation
where $\CP^2$ is blown up at nine points which lie at the intersection of
two elliptic curves. In this case the mirror is an elliptic fibration
over $\C$ for which the compactification has a smooth fiber at infinity.
Compared to that of $\CP^2$ ($k=0$), this extreme case where $k=9$ lies
at the opposite end of the spectrum that we consider.
\end{remark}


\end{document}